\documentclass{article}

\usepackage[a4paper,margin=1in,footskip=0.25in]{geometry}

\usepackage[utf8]{inputenc} 
\usepackage[T1]{fontenc}    
\usepackage{hyperref}       
\usepackage{url}            
\usepackage{booktabs}       
\usepackage{amsfonts}       
\usepackage{nicefrac}       
\usepackage{microtype}      
\usepackage{csquotes}
\usepackage{amsthm}
\usepackage{subfig}

\usepackage{graphicx}
\usepackage{amsmath,amssymb}
\usepackage{tikz}

\usetikzlibrary{arrows}
\usetikzlibrary{calc,arrows.meta,positioning}

\def\bydef{\triangleq}

\def\N{\mathbb{N}}
\def\R{\mathbb{R}}
\def\rp{\mathbb{R}_+}

\def\calA{\mathcal{A}}

\def\calI{\mathcal{I}}

\def\calY{\mathcal{Y}}

\def\calS{\mathcal{S}}
\def\calP{\mathcal{P}}

\def\sol{{\rm ALG}}
\def\tc{\pi_{\rm TC}}
\def\tv{\delta_{\rm TV}}

\def\E{\mathbb{E}}

\def\nln{\nonumber\\}

\newtheorem{definition}{Definition}
\newtheorem{theorem}{Theorem}
\newtheorem{lemma}{Lemma}

\makeatletter
\tikzset{my loop/.style =  {to path={
  \pgfextra{}
  [looseness=6,min distance=1mm]
  \tikz@to@curve@path},font=\sffamily\small
  }}  
\makeatletter

\title{Temporal Concatenation for Markov Decision Processes}

\author{%
  Ruiyang Song \\
  Department of Electrical Engineering\\
  Stanford University\\
  \texttt{ruiyangs@stanford.edu} \\
  \and
   Kuang Xu\\
   Graduate School of Business\\
  Stanford University \\
   \texttt{kuangxu@stanford.edu} \\
}

\begin{document}

\maketitle
\pagenumbering{arabic}
\begin{abstract}

We propose and analyze a temporal concatenation heuristic for solving large-scale finite-horizon Markov decision processes (MDP), which divides the MDP into smaller sub-problems along the time horizon and generates an overall solution by simply concatenating the optimal solutions from these sub-problems. As a “black box” architecture, temporal concatenation works with a wide range of existing MDP algorithms. Our main results characterize the regret of temporal concatenation compared to the optimal solution. We provide upper bounds for general MDP instances, as well as a family of MDP instances in which the upper bounds are shown to be tight. Together, our results demonstrate temporal concatenation's potential of substantial speed-up at the expense of some performance degradation.

\end{abstract}

\section{Introduction}
We are interested in devising computationally efficient architectures for solving finite-horizon Markov decision processes (MDP), a popular framework for modeling multi-stage decision-making problems \cite{bellman1954theory, puterman2014markov} with a wide range of applications from scheduling in data and call centers \cite{harrison2004dynamic} to energy management with intermittent renewable resources \cite{harsha2015optimal}. In an MDP, at each stage, an agent makes a decision based on the state of the system, which leads to an instantaneous reward and the state is updated accordingly; the agent aims to find an optimal policy that maximizes the total expected rewards over the time horizon.  While finding efficient algorithms for solving MDPs has long been an active area of research (see~\cite{mundhenk2000complexity, littman1995complexity} for a survey), we will, however, take a different approach. Instead of creating new algorithms from scratch, we ask how to design \emph{architectures} that leverage existing MDP algorithms as ``black boxes''  in creative ways, in order to harness additional performance gains. 

As a first step in this direction, we propose the temporal concatenation heuristic, which takes a divide-and-conquer approach along the time axis: for an MDP with horizon $\{0, \ldots, T-1\}$, we divide the original problem instance ($\calI_0$) over the horizon into two sub-instances:  $\left\{0, \ldots, \frac{T}{2}-1\right\}$ ($\calI_1$) and $\left\{\frac{T}{2}, \ldots, T-1\right\}$ ($\calI_2$), respectively. Temporal concatenation then evokes an MDP algorithm, one that takes as input an MDP instance and outputs an optimal policy, to find the optimal policies $\pi_1^*$ and $\pi_2^*$ for the two sub-instances $\mathcal I_1$ and $\mathcal I_2$, separately. Finally, temporal concatenation outputs a policy, $\tc$, for the original MDP by simply concatenating $\pi^*_1$ and $\pi^*_2$: run $\pi^*_1$ during the first half of the horizon, and $\pi^*_2$ the second.\footnote{ More generally, a similar temporal concatenation procedure can be performed over $K$ sub-instances, with $K \geq 2$. Our theoretical analysis will focus on the case of $K=2$ because it captures the majority of structural insights. 
}

In a nutshell, temporal concatenation is intended as a simple ``black box'' architecture to substantially speed up existing MDP algorithms, at the expense of potentially minor performance degradation.
First, acceleration comes from the fact that the optimal policies for the sub-instances can be derived entirely in parallel.
In particular, a classical MDP problem can be solved by conventional methods such as the value iteration, which has a time complexity growing linearly with the horizon, $T$. By applying the temporal concatenation architecture in this set-up, the computing time can, in principle, be reduced by half.  
 This speed-up from parallelism can be significant if the original MDP algorithm's run-time is sufficiently long. 
 In addition, we point out that a parallelism architecture is well suited for modern machine learning systems where the instance of a large-scale problem may be stored in separate servers to start with \cite{low2012distributed,chen2015mxnet, meng2016mllib}. 
 Moreover, temporal concatenation can speed up computation even more significantly if the MDP algorithm in question admits a run-time that scales super-linearly in the horizon. Typically, these algorithms suffer a worse dependence on $T$ in exchange for a more favorable scaling in the size of the state and action spaces; see for instance, the complexity of a linear-programming-based algorithm  in \cite{tseng1990solving}  that scales as $\mathcal{O}(T^4)$, and that of the stochastic primal-dual method proposed by Chen and Wang \cite{chen2016stochastic}, which scales as $\mathcal{O}(T^6)$.

While the computational benefit from using temporal concatenation is evident, the quality of its solution is not: by solving two sub-instances independently, it could be overly short-sighted and lead to strictly sub-optimal MDP policies. Therefore, our theoretical results will focus on addressing the following question: 

\emph{How good is the policy generated by temporal concatenation, $\pi_{\rm TC}$, compared to the optimal policy to the original problem, $\pi^*$? }

{\bf Preview of main results}. On the positive side, we provide sufficient conditions under which the performance gap between $\pi_{\rm TC}$ and $\pi^*$ is \emph{small}. Specifically, we establish upper bounds to show that the performance gap is bounded by a function that depends linearly on an MDP's \emph{diameter} (a measure that reflects the ease with which the agent can traverse the state space) but \emph{independent} from the horizon, $T$. Conversely, we provide lower bounds by showing that, for \emph{any} finite diameter, there exist MDP instances for which the upper bounds are tight for all large $T$.

{\bf Organization}. The remainder of this paper is organized as follows. In Section \ref{sec:problem_formulation}, we formally introduce the problem formulation and performance metrics. In Section \ref{sec:main_result}, we summarize the main results and contrast our approach to the extant literature. 
Section \ref{sec:examples} provides several examples of MDP instances, including one that is motivated by the application of dynamic energy management with on-site storage. We also provide simulation results that substantiate the theoretical results and illustrate the run-time reduction obtained by running the temporal concatenation heuristic on a multi-core PC. Section \ref{sec:conclusion} concludes the paper.

{\bf Notation}. We will denote by $[n]$ the set of integers $\{0,1,\ldots,n-1\}$, $n\in\mathbb{N}$. We will use $\delta_{\rm TV}(\mu,\nu)$ to denote the total variation distance between two distributions $\mu$ and $\nu$: $\delta_{\rm TV}(\mu,\nu) = \frac{1}{2}\sum_{s}|\mu(s)-\nu(s)| = \sum_{s: \mu(s) \geq \nu(s)}\left(\mu(s)-\nu(s)\right)$. For a sequence $\{a_i\}_{i\in \N}$, and $s,t\in \N$, $s\leq t$,  we use $s\to t $ to denote the set $\{s,s+1,\ldots,t\}$, and use $a_{s\to t}$ to denote the sub-sequence $\{a_s,a_{s+1}, \ldots, a_{t-1},a_t\}$. Similarly, for some $\calS\subseteq \N$, $a_\calS$ denotes the set $\{a_i: i\in \calS\}$. 
For $x\in \R$, we will denote by $(x)^+$ and $(x)^-$ the positive and negative portion of $x$, respectively: $(x)^+ = \max\{x,0\}$ and $(x)^- = \max\{-x,0\}$. For $c,d\in \mathbb{R}$ with $c\leq d$, define $x_{[c,d]}$ to be the projection of $x$ onto the interval $[c,d]$, i.e., $x_{[c,d]}\triangleq \mathbb{I}(x< c)c+\mathbb{I}(c\leq x\leq d)x+\mathbb{I}(x> d)d$, where $\mathbb{I}(\cdot)$ is the indicator function.

\section{Problem Formulation and Performance Metric}
\label{sec:problem_formulation}

\subsection{System Set-up} 
We consider a discrete-time Markov decision process with a finite time horizon $[T]$, state space $\calS$, and action set $\calA$. The decision maker chooses at each step $t \in [T]$ an action, $a_t \in \calA$.  We will assume that $\calA$ and $\calS$ stay fixed, and hence omit them from our notation when appropriate. The state of the system at time $t$ is denoted by $S_t$. The initial state $S_0$ is drawn from some probability distribution $\mu_0$, and the state evolution depends on the present state as well as the action chosen: 
\begin{equation}
S_{t+1}=p_t\left(a_t,S_t,Y_t^S\right),\quad t\in[T].
\end{equation}
The $Y^S_t$'s are i.i.d.~uniform random variables over a finite set $\mathcal Y^S$, capturing the randomness in the state transition. The collection $\{p_t\}_{t\in [T]}$ is the set of (deterministic) transition functions. The decision maker receives a \emph{reward} at each time slot $t$,  $R_t(a_t,S_t, Y^R_t)$, which depends on the present state, action, and some i.i.d.~idiosyncratic random variables taking values in a finite set, $Y^R_t\in\mathcal Y^R$, with a fixed distribution. We refer to $\{R_t\}_{t\in [T]}$ as the set of reward functions. We assume that the rewards are nonnegative and bounded from above by a constant, $\bar{r} \in \rp$. 
        \footnote{{Note that the results in this paper would be unchanged if the reward function $R_t$ were shifted by a constant. In particular, the general case in which $R_t(\cdot,\cdot,\cdot)\in[r_{\min},r_{\max}]$ for $-\infty<r_{\min}\leq r_{\max}<\infty$, is equivalent to having $R_t\in[0,\bar{r}]$ where $\bar{r}=r_{\max}-r_{\min}$. Throughout this paper we let $R_t(\cdot,\cdot,\cdot)\in[0,\bar{r}]$ for simplicity of notation unless otherwise specified.}} 
        The decision maker's  behavior is described by a \emph{policy} $\pi(\cdot)$, such that $a_t = \pi(t,S_t, Y^P)$, $ t\in [T]$. In other words, the policy chooses an action based on the current state, and some idiosyncratic randomization $Y^P$, which, without loss of generality, can be thought of as a uniform random variable over $[0,1]$.

An MDP as described above is specified by the triple, $(R_{[T]}, p_{[T]},T)$, which we will refer to as a \emph{problem instance}. We will refer to the original horizon-$T$ MDP problem instance as the \emph{original instance}, denoted by $\calI_0 \bydef (R_{[T]}, p_{[T]},T). $ For an instance $\calI_0$, policy $\pi$, and initial distribution $\mu_0$, the \emph{total expected reward} is defined by
\begin{equation}
V(\calI_0, {\pi}, \mu_0)= \E^{\pi}_{S_0 \sim \mu_0}\left[\sum_{t=0}^{T-1} R_t\left(a_t,S_t,Y^R_t\right)\right], 
\label{eq:defReg}
\end{equation}
where the expectation is taken with respect to all the randomness in the system, and the actions are chosen according to $\pi$. A policy $\pi$ is \emph{optimal} if it attains the maximum total expected reward for all initial distributions, $\mu_0$.

\subsection{Temporal Concatenation} 
We now define the main object of study, the temporal concatenation heuristic. An \emph{MDP algorithm}, denoted by $\sol(\cdot)$, takes as input a problem instance and outputs the optimal policy, $\pi^*$, for that instance. As such, the notion of an MDP algorithm captures the ``functionality'' of an algorithm that is used to compute an optimal policy,  but abstracts away the inner working of the algorithm, effectively treating it as a ``black box.'' By this definition, we have that
\begin{equation}\label{eq:solver}
\pi^* = \sol(\calI_0). 
\end{equation}
\begin{definition}[Temporal concatenation] \label{def:tpc}
For an {\bf original instance} $\calI_0$, denote by $\calI_1$ and $\calI_2$ the {\bf sub-instances} generated by partitioning  $\calI_0$ in half along the time horizon: 
\begin{equation}\label{def:sub-instances}
\calI_1 \bydef \left(R_{0\to \frac{T}{2}-1}, p_{0\to  \frac{T}{2}-1}, \frac{T}{2}\right), \quad \mbox{and} \quad \calI_2 \bydef \left(R_{\frac{T}{2} \to T-1}, p_{\frac{T}{2} \to T-1}, \frac{T}{2}\right),
\end{equation}
and by $\pi^*_1$ and $\pi^*_2$ their corresponding optimal policies: 
\begin{equation}
\pi^*_1\bydef \sol(\calI_1), \quad \mbox{and} \quad \pi^*_2\bydef \sol(\calI_2). 
\end{equation} 
The temporal concatenation heuristic generates a policy, $\tc$, by temporally concatenating optimal solutions for $\calI_1$ and $\calI_2$, $\pi^*_1$ and $\pi^*_2$, i.e., 
\begin{equation}\label{eq:tpc1}
\tc\left(t, S_t, Y^P\right) = \left\{
                \begin{array}{ll}
                  \pi^*_1\left(t,S_t, Y^{P}_1\right), & \quad 0\leq t\leq \frac{T}{2}-1, \\
                   \pi^*_2\left(t-T/2,S_t, Y^{P}_2\right), & \quad \frac{T}{2} \leq t\leq T-1,
                \end{array}
              \right.
\end{equation}
where $Y_1^P$ and $Y_2^P$ are two independent uniform random variables generated from $Y^P$ as follows: express $Y^P \in [0,1]$  as an infinite binary sequence, and set $Y_1^P$ and $Y_2^P$ to be the sub-sequence corresponding to all odd and even elements in the binary sequence, respectively.
\end{definition}

\subsection{Performance Metric} 
The following definition of regret is our main metric, which measures how the expected reward of the policy $\tc$ deviates from the optimal policy $\pi^*$: 
\begin{definition}[Regret of temporal concatenation]\label{def:regret}
For an original instance $\calI_0$ and initial distribution $\mu_0$, the \emph{regret of temporal concatenation}, or regret for short,  is defined by: 
\begin{equation}
\Delta(\calI_0, \mu_0) \bydef  V({\calI_0}, \pi^*, \mu_0) - V({\calI_0}, \tc, \mu_0),  
\end{equation}
where $\pi^*$ is an optimal policy for $\calI_0$, and $\tc$ is defined in Definition \ref{def:tpc}.  
\end{definition}
Note that the above definition differs from the conventional notion of regret of the online learning literature and reinforcement learning. For instance, regret in \cite{burnetas1997optimal} is due to not having complete information of the MDP in hindsight, whereas in our case, it is due to the intrinsic sub-optimality from dividing an original MDP instance into smaller sub-problems and solving each separately.

\section{Main Results}\label{sec:main_result}

We present our main results in this section. The first result, Theorem \ref{thm:upperbound}, provides an upper bound on the regret of temporal concatenation in an MDP, which does not depend on the length of the horizon, $T$. Instead, the regret is shown to be related to a notion of \emph{diameter} of the MDP, which we define below.  

The diameter captures how easy it is for the decision maker to reach different state distributions. Let $\mathcal P$ be the collection of all distributions over $\mathcal S$. 
We will denote by $\mu_t^\pi$ the state distribution at time $t$ under policy $\pi$.
Starting at time $t_0\geq 0$, for two distributions $\mu, \nu\in\mathcal P$, we say that $\nu$ is \emph{$\epsilon$-reachable} from $\mu$ in $t$ steps for some $\epsilon\in[0,1]$, if there exists a policy $\pi$ such that under $\pi$ and with the distribution of $S_t$ at time $t_0$ being $\mu$, we have that 
\begin{equation}
\tv(\mu_{t_0+t}^\pi, \nu) \leq \epsilon.
\end{equation}
 Denote by $\calP_\epsilon^{t_0}(\mu,t)$ the set of all distributions that are $\epsilon$-reachable from $\mu$ in $t$ steps starting from time $t_0\in [T-t]$.  We have the following definition of diameter. 

\begin{definition}[$\epsilon$-Diameter]\label{def:eps-diameter}
For an MDP instance $\calI_0$ with horizon $[T]$ and transition functions $\{p_t\}_{t\in [T]}$, we define the $\epsilon$-diameter as the least number of steps with which, starting from any time step, all possible distributions in $\mathcal P$ are $\epsilon$-reachable from one another:
\begin{align}
&\tau_\epsilon(\calI_0) \triangleq  \inf\left\{t\geq 0: \nu'\in\mathcal P^{t_0}_\epsilon(\nu,t)\ {\rm for\ all\ } \nu,\nu'\in\mathcal P{\rm \ and\ all\ }t_0\in[T-t]\right\}.
\end{align}
\end{definition}
Note that since $\tau_\epsilon$ captures the hardness of traversing the state space by applying feasible actions, it will depend on the sizes of the state  and action spaces in general.
We have the following theorem. 
\begin{theorem}[Upper bound on regret of temporal concatenation]
\label{thm:upperbound}
Fix an original instance $\calI_0$ with horizon $[T]$, and an initial distribution $\mu_0$.  If there exists $\epsilon\geq 0$ such that $\tau_\epsilon(\calI_0)\leq T/2$,  then the regret of temporal concatenation satisfies: 
\begin{equation}
\Delta(\calI_0, \mu_0)\leq \frac{\bar{r}\tau_\epsilon(\mathcal I_0)}{1-\epsilon},
\label{eq:thm_upper_bound}
\end{equation}
where $\bar{r}$ is the maximum reward in a given time slot.  In particular, if $\tau_0(\calI_0)\leq T/2$,  then the above inequality implies that
\begin{equation}\label{eq:tau_0}
\Delta(\calI_0, \mu_0) \leq \bar{r} \tau_0(\calI_0). 
\end{equation}
\end{theorem}

A direct implication of the above theorem is that, for problems that admit a moderate $\epsilon$-diameter for some $\epsilon\in [0,1)$, temporal concatenation produces a near-optimal policy \emph{regardless} of the length of the horizon, $T$, thus making the heuristic especially appealing for problems with a relatively large horizon. 

It is also worth noting that while the original temporal concatenation algorithm requires an optimal policy to be used in each of the two sub-instances, it is easy to substitute these optimal policies with sub-optimal ones with a bounded regret and obtain similar regret bounds to those in Theorem \ref{thm:upperbound}. In particular, suppose we use in each sub-instance a policy whose total expected reward over the sub-instance is at most $\delta$ less than that of an optimal policy starting from any initial state, then it follows the resulting regret bounds would  be those in Theorem \ref{thm:upperbound} with an additional $ 2\delta $ additive factor. 

The next result provides a lower bound that demonstrates that a small diameter is also \emph{necessary} for temporal concatenation to perform well,  in a worst-case sense. We look at MDP instances with a bounded $0$-diameter, $\tau_0(\mathcal I_0)$. In Theorem \ref{thm:lowerb_finite_diam}, we show that, for any $d_0\in \N$, there exists an instance with a $0$-diameter equal to $d_0$ such that the performance regret is essentially $\bar{r} d_0 $ for any horizon $T>2d_0+2$. This result implies that the upper bound in \eqref{eq:tau_0} of Theorem \ref{thm:upperbound} is tight in a worst-case instance.   
\begin{theorem}[Lower bound on regret of temporal concatenation]
\label{thm:lowerb_finite_diam}
Fix $\bar{r}\in\mathbb R_+$, $\sigma\in(0,\frac{\bar{r}}{2})$, and integer $d_0\geq 5$. Then there exists an MDP instance $\calI_0$ with maximum per-slot reward $\bar{r}$, finite $0$-diameter $\tau_0(\calI_0)=d_0$, and an initial distribution $\mu_0$, such that for any $T> 2d_0+2$, the regret satisfies   
\begin{equation}
\Delta(\mathcal I_0,\mu_0) = (\tau_0(\calI_0)-2)\bar{r}-\sigma.
\end{equation}
\end{theorem}

 The proofs of Theorems \ref{thm:upperbound} and \ref{thm:lowerb_finite_diam} will be presented in Appendix \ref{sec:proof}.

Theorem \ref{thm:upperbound} shows that for MDP instances that admit a bounded $\epsilon$-diameter, the  regret of temporal concatenation is bounded from above by a value independent of the time horizon $T$. This is encouraging, since it would suggest that the quality of approximation afforded by our heuristic does not degrade over longer time horizons. On the other hand, Theorem \ref{thm:lowerb_finite_diam} shows that the regret could be very large if the diameter is large, though it would appear from our proof that the ``bad'' examples we know so far tend to be fairly pathological. 
We accompany these findings  by examining in Section \ref{sec:examples} a number of specific MDP models. Our theoretical and simulation results there suggest that temporal concatenation at least  performs reasonably well for several such ``average'' instances. 

\subsection{Related Work}

Our method is related to the literature on MDP decomposition methods, which aim to  overcome the so-called curse of dimensionality by breaking down the original MDP with a large state space into sub-problems with smaller state spaces. Hierarchical MDP algorithms utilize hierarchical structures to decompose the state space and action space and transform the original problem into a collection of sequential sub-tasks \cite{parr1998reinforcement, dietterich2000hierarchical, mundhenk2000complexity, daoui2010exact}. The method in \cite{sucar2007parallel} decomposes the problem into parallel sub-tasks that can be computed simultaneously, where each sub-task is an MDP with the same state and action space, but with different reward functions. Steimle \textit{et al.}~\cite{steimle2019decomposition} adopts a decomposition via mixture models. Finally, Ie \textit{et al.}~\cite{ie2019slateq} leverages decomposition in the Q-function by exploiting combinatorial structures of the recommender systems.

While our approach also works by dividing the original MDP instance into smaller sub-problems,  there is a number of crucial features that differentiate our approach.   Firstly, our method focuses on decomposing the problem along the time axis rather than over the state space or action space, which is a more common approach in the literature.  Secondly, the focus of temporal concatenation  is to serve as a simple ``black box'' architecture, rather than a custom-made MDP algorithm.  As such, each sub-problem can be solved by any MDP algorithm of the user's choice, and the procedure is very simple to implement and does not involve complex procedures to transform the structure of the original problem. Finally, as alluded to in the introduction, temporal concatenation lends itself easily to parallel processing, and thus achieving speed-up that is not possible under a decomposition algorithm in which sub-problems still need to be solved in a sequential manner (cf.~\cite{sucar2007parallel}).

Related to our approach in spirit is \cite{kolobov2012reverse}, which proposes a heuristic for finite-horizon MDPs by sequentially solving a series of smaller MDPs with increasing horizons, and the numerical results show that the heuristic provides good performance even if the process is terminated prematurely. However, no rigorous guarantees in terms of regret of this heuristic relative to the optimal policy were established.

Related notions of the diameter of an MDP have been used in the literature to capture the ease with which the system can transit between any pair of states in the state space. For instance, a diameter $D^*$ is defined in \cite[Definition 1]{jaksch2010near} as  
\begin{equation}
D^* \triangleq \max_{s,s'\in\mathcal S,\ s\neq s'} \min_{\pi}\E^{\pi}\left[\min_{N\geq 1,S_N = s}N\Big | S_0 = s\right].
\end{equation}
 The diameter $D^*$ has been used for analyzing the total regret of reinforcement learning algorithms (see \cite{jaksch2010near,talebi2018variance} for example). In \cite{jaksch2010near}, the authors introduced a learning algorithm for MDP with total regret $\mathcal O(D^*|\mathcal S|\sqrt{|\mathcal A|T})$. In \cite{talebi2018variance}, an improved upper bound for the total regret of MDP is introduced, which depends on the variance of the bias function defined in \cite[Definition 2]{talebi2018variance}, but does not depend on $D^*$.
In comparison, our definition of $\epsilon$-diameter is different, and in some sense stronger. 
First, while $D^*$ is the expected number of steps necessary to transit between any pair of states, $\epsilon$-diameter corresponds to the number of steps required to traverse between any pair of \emph{distributions} over the state space. Second, the definition of $D^*$ implies the existence of a policy under which the target state is reached with no more than $D^*$ time steps, while for $\tau_\epsilon$, we require a policy such that the target distribution is achieved after exactly $\tau_\epsilon$ time steps (a total variation distance no greater than $\epsilon$ is allowed). Notably, the lower bound in Theorem \ref{thm:lowerb_finite_diam} shows that our notion of diameter cannot be weakened when applied to the analysis of temporal concatenation, thus suggesting that our formulation reveals structural properties of the MDP distinct from those in the extant literature. 
We will further discuss the connection between the $\epsilon$-diameter and the $D^*$ diameter in Appendix \ref{sec:connection_to_d_star}.

\section{Examples and Illustrative Applications}
\label{sec:examples}

In this section we discuss several examples to illustrate the properties of the $\epsilon$-diameter and corroborate the theoretical results in Section \ref{sec:main_result}. In Section \ref{subsec:dgt}, we introduce the deterministic graph traversal (DGT) problems, a family of MDP instances with finite $0$-diameter and noiseless transitions. In Section \ref{subsec:xisgt}, we introduce the $\xi$-stochastic graph traversal ($\xi$-SGT) problems, which is a generalization of the DGT with stochastic transitions. In Section \ref{subsec:dem}, we present a model of dynamic energy management with storage, which is an illustrative example of the $\xi$-SGT family. We also present simulation results for the deterministic graph traversal models in Section \ref{subsec:numerical} to explore the average-case scaling behavior of the regret within this family. In Section \ref{subsec:runtime}, we provide additional simulation examples of this model to illustrate the run-time reduction from using the temporal concatenation heuristic. In Section \ref{subsec:garnet}, we present simulation results of a more popular family of MDP instances known as the Generalized Average-Reward Non-stationary Environment Test-bench (GARNET) model.

\subsection{Deterministic Graph Traversal Problems}\label{subsec:dgt}
In this subsection, we introduce the \emph{deterministic graph traversal (DGT)} problems, a family of MDP instances with finite $0$-diameter. Let $\mathcal G_{\rm csl}$ be the set of all strongly connected graphs that include at least one self-loop.
A DGT instance denoted by $\mathcal I_G$ has a time-homogeneous deterministic transition function, i.e. $p_t=p$ for all $t$, which can be described by a strongly connected directed graph $G=(\mathcal V,\mathcal E)\in\mathcal G_{\rm csl}$ where at least one vertex in $G$ has a self-loop. Here, $\mathcal V$, $\mathcal E$ are the collections of vertices and edges of $G$, respectively. 
In other words, for any $G=(\mathcal V,\mathcal E)\in\mathcal G_{\rm csl}$, there exists a vertex $v_i\in\mathcal V$ such that the self-loop edge exists, i.e., $e_{ii}\in\mathcal E$.
In a DGT instance, once the current state $S_t$ and action $a_t$ are given, the next state $S_{t+1}$ is determined. 
We formally define a DGT instance with state space $\mathcal S$, action space $\mathcal A$, and transition function $p$ as follows.
\begin{definition}[Deterministic graph traversal instance]\label{def:dgt}
Let $G=(\mathcal V,\mathcal E)\in\mathcal G_{\rm csl}$.
 A DGT instance depicted by $G$, $\mathcal I_G$, is an MDP instance whose state space and transition function satisfy:
         
         (1) Each state $i\in\mathcal S$ corresponds to a vertex $v_i\in\mathcal V$. 
         
         (2) The state transition is deterministic and can be described by the edges in $\mathcal E$.
         In particular, for $i,j\in \mathcal S$, an edge $e_{ij}\in\mathcal E$ implies the existence of an action $a_{ij}\in\mathcal A$ such that  starting from state $i$, the system will deterministically go to state $j$ once the agent takes action $a_{ij}$, i.e., $j = p(a_{ij},i,Y_t^S)$ with probability 1 for all $t$. 
\end{definition}
Note that we are not imposing additional restrictions on the reward functions in the definition of DGT instances. 
As an example, the MDP instance we construct for proving Theorem \ref{thm:lowerb_finite_diam} is a special case that belongs to this family (see Appendix \ref{sec:proof_lower_fin_Diam}).

Now we study the $\epsilon$-diameter of DGT instances. 
We first briefly recall the definition of the classical diameter of a directed graph $G= (\mathcal V, \mathcal E)$, denoted by $d_c(G)$. 
For any two vertices of the graph, $i,j\in\mathcal V$, let $d_G(i,j)$ be the distance between them on graph $G$, which is defined as the length of the shortest path from $i$ to $j$. Here, a path is a sequence of distinct vertices such that each two consecutive vertices are connected by an edge in $\mathcal E$. The classical diameter is the maximum taken over all pairwise distances, i.e., $d_c(G) = \max_{i,j\in\mathcal V}d_G(i,j)$. 
For a strongly connected graph $G$, each pairwise distance is finite, in which case $d_c(G)$ is also finite. Further, it is not difficult to verify that for a strongly connected graph, the classical diameter is at most $|\mathcal V|-1$. Literature has shown that $d_c(G)$ can be computed within at most $O(|\mathcal V|^3)$ time using classical algorithms such as the breadth first search (see \cite{bondy1976graph} for example).

Recall that $d_c(G)<\infty$ because $G$ is strongly connected. In the following lemma, we prove that DGT instances indeed have a finite $0$-diameter, which is closely related to the classical diameter of the corresponding graph, $d_c(G)$. 
\begin{lemma}[From $\epsilon$-diameter to classical graph diameter] \label{lemma:connection of diameters}
For a DGT instance, $\mathcal I_G$, based on a graph $G=(\mathcal V,\mathcal E)\in\mathcal G_{\rm csl}$, we have 
\begin{equation}
d_c(G)\leq \tau_0(\mathcal I_G)\leq 2d_c(G).
\end{equation}
\end{lemma}
With Lemma \ref{lemma:connection of diameters}, we have shown that the $0$-diameter of $\mathcal I_G$ is finite and bounded between $d_c(G)$ and $2d_c(G)$. It follows that for any $\epsilon>0$, the $\epsilon$-diameter of $\mathcal I_G$ satisfies
\begin{equation}
\tau_\epsilon(\mathcal I_G)\leq \tau_0(\mathcal I_G)\leq 2d_c(G).
\end{equation}
The proof of this lemma is given in Appendix \ref{app:lemma:connection of diameters}.

\subsection{$\xi$-Stochastic Graph Traversal Problems}\label{subsec:xisgt}
In Section \ref{subsec:dgt}, we introduced a family of MDP instances with finite $0$-diameter $\tau_0$ where the transition is deterministic.
However, this model can not capture the stochasticity in many real-world applications. 
To this end, we study in this subsection a generalization of the DGT family where transitions can be impacted by stochastic shocks. As a result, we will see concrete examples of instances where the $\epsilon$-diameter is finite for $\epsilon>0$, even though the $0$-diameter may be infinite.  

Specifically, we consider the \emph{$\xi$-stochastic graph traversal ($\xi$-SGT)} problems with state space $\mathcal S$, action space $\mathcal A$, and time-homogeneous transition function $p$, defined as follows.
\begin{definition}[$\xi$-Stochastic graph traversal instance]\label{def:xisgt}
Fix $G=(\mathcal V,\mathcal E)\in\mathcal G_{\rm csl}$, and $\xi\in(0,1)$. For $i\in\mathcal V$,  define the neighbor set of $i$ as $\mathcal N_i =\{j\in\mathcal V: (i,j)\in \mathcal E \}$. A $\xi$-SGT instance based on $(G,\xi)$, $\mathcal I_G^\xi$, is an MDP instance whose state space and transition function satisfy:
         
         (1) Each state $i\in\mathcal S$ corresponds to a vertex $v_i\in\mathcal V$. 
         
         (2) The state transition is stochastic and can be described by the edges in $\mathcal E$. In particular, for $i,j\in \mathcal S$ such that the edge $e_{ij}\in\mathcal E$ exists, there exists an action $a_{ij}\in\mathcal A$ under which starting from state $i$, the system goes to state $j$ with a probability at least $1-\xi$, i.e., 
         \begin{equation}
         p\left(a_{ij},i,Y_t^S\right)=
         \begin{cases}
         j & {\rm with\ a \  probability\ at\ least}\ 1-\xi \\
         Z_{t,i,j} & {\rm otherwise}
         \end{cases},
          \end{equation}
          where $Z_{t,i,j}$ is a random variable that takes values in $\mathcal N_i\backslash\{j\}$ if $|\mathcal N_i|\geq 2$, and $Z_{t,i,j}  = j$ if $|\mathcal N_i| = 1$.
\end{definition} 
From Definitions \ref{def:dgt} and \ref{def:xisgt}, we see that the $\xi$-SGT and DGT problems are closely related. For both of these families of MDP instances, the state space corresponds to the vertices of a directed graph, and the state transition function can be described by the edges in the same graph. 
The $\xi$-SGT instance can be regarded as a noisy version of the DGT instance, where the transition along the edges is perturbed by a random noise. The parameter $\xi$ can be interpreted as the noise level. 
In a DGT instance based on $G$, the system can deterministically traverse the state space along the edges of $G$ when appropriate actions are taken. In an $\xi$-SGT based on $(G,\xi)$, however, when a proper action is chosen, the system will traverse along the ``intended'' edge with a probability at least $1-\xi$, but may be diverged to one of the other neighbors otherwise. 
The following result connects the $\epsilon$-diameter of the $\xi$-SGT instance, $\tau_\epsilon\left(\mathcal I_G^\xi\right)$, the $0$-diameter of the DGT instance, $\tau_0\left(\mathcal I_G\right)$, and the classical diameter of the underlying graph, $d_c(G)$; the proof is given in Appendix \ref{app:lemma:xisgt diameter}.
\begin{lemma}[Diameter of $\xi$-SGT instances]\label{lemma:xisgt diameter}
Fix $G\in \mathcal G_{\rm csl}$, and $\xi\in(0,1)$. Let $\mathcal I_G$ be the DGT instance characterized by $G$, and $\mathcal I_G^\xi$ the $\xi$-SGT instance described by $(G,\xi)$. For $\epsilon\geq 1-(1-\xi)^{\tau_0(\mathcal I_G)}$, we have
\begin{equation}
\tau_\epsilon\left(\mathcal I_G^\xi\right) \leq\tau_0(\mathcal I_G).
\end{equation}
Combining the above with Lemma \ref{lemma:connection of diameters}, we have 
\begin{equation}\label{eq:eps-classical}
\tau_\epsilon\left(\mathcal I_G^\xi\right) \leq 2d_c(G)
\end{equation}
for $\epsilon\geq 1-(1-\xi)^{2d_c(G)}$.
\end{lemma}
 Lemma \ref{lemma:xisgt diameter} implies the following communicating property of a $\xi$-SGT instance: when we are allowed a total variation distance $\epsilon$ and the noise level $\xi$ is sufficiently small such that $\xi<1-(1-\epsilon)^{\frac{1}{\tau_0(\mathcal I_G)}}$, we can traverse the state space in a $\xi$-SGT instance using no more than $\tau_0(\mathcal I_G)$ steps. 
Further, with Lemma \ref{lemma:connection of diameters}, the $0$-diameter of the DGT instance, $\mathcal I_G$, is bounded from above by two times the classical diameter of the corresponding graph $G$. This in turn implies that the $\epsilon$-diameter is bounded from above by $2d_c(G)$.

While determining the closed-form $0$-diameter of a general DGT instance remains a direction for future work, we have provided in this paper an example of DGT instance whose $\tau_0$ can be derived explicitly. In particular, the MDP instance depicted by Figure~\ref{fig:prop_large_regret} in Appendix~\ref{sec:proof_lower_fin_Diam}, which is designed for proving the lower bound in Theorem \ref{thm:lowerb_finite_diam}, is a DGT instance of $k+2$ states with    
$\tau_0(\mathcal I_G)= k+2$. 
For $\xi$-SGT instances corresponding to the graph in Figure \ref{fig:prop_large_regret}, we can apply Lemma \ref{lemma:xisgt diameter} and obtain $\tau_\epsilon(\mathcal I_G^\xi)\leq k+2$ for $\epsilon \geq 1-(1-\xi)^{k+2}$. 
If the 0-diameter of the DGT instance cannot be derived in closed form, we can apply the upper bound in Eq.~\eqref{eq:eps-classical} instead. This upper bound depends only on the classical diameter, which can be calculated for any directed graph by existing algorithms.

 While Lemma \ref{lemma:xisgt diameter} is a general result that holds for any $\xi$-SGT instance, we introduce another stronger characterization on the $\epsilon$-diameter of $\xi$-SGT instances when the graph $G=(\mathcal V,\mathcal E)\in \mathcal  G_{\rm csl}$ is undirected, and each node of $G$ has a noiseless self-loop, i.e., $(i,i)\in\mathcal E$, and $p(a_{ii},i,Y_t^S)=i$ with probability one, for all $i\in\mathcal S$. We have the following lemma, which is proved in Appendix \ref{app:lemma:xisgt new}.
 
\begin{lemma}\label{lemma:xisgt new}
Fix an undirected connected graph $G=(\mathcal V,\mathcal E)\in\mathcal G_{csl}$ where each node has a noiseless self-loop, i.e., $(i,i)\in\mathcal E$ for all $i\in\mathcal S$, and $\xi\in\left(0,\frac{1}{2}\right)$. Let $\mathcal I_G$ be the DGT instance characterized by $G$, and $\mathcal I_G^\xi$ be the $\xi$-SGT instance described by $(G,\xi)$. Then for any $\epsilon\in(0,1)$, the $\epsilon$-diameter of $\mathcal I_G^\xi$ satisfies that 
\begin{equation}\label{eq:lemma xisgt new}
\tau_\epsilon\left(\mathcal I_G^\xi\right)\leq \frac{d_c(G)}{1-2\xi}+\frac{f(d_c(G),\xi)}{\epsilon},
\end{equation}
with 
\begin{equation}\label{eq:def:f}
f(d_c(G),\xi) = \frac{4\xi(1-\xi)}{(1-2\xi)^2}\left(2+\sqrt{4+\frac{(1-2\xi)d_c(G)}{\xi(1-\xi)}}\right),
\end{equation} 
where $d_c(G)$ is the classical diameter of graph $G$.
Note that when $\xi\to 0$, Eq.~\eqref{eq:lemma xisgt new} reduces to $\tau_\epsilon\left(\mathcal I_G\right)\leq d_c(G)$.
\end{lemma}
 Lemma \ref{lemma:xisgt new} suggests that for a $\xi$-SGT instance with undirected $G$ and noiseless self-loops for all nodes, there is an upper bound on the $\epsilon$-diameter, which grows linearly in $\frac{1}{\epsilon}$. Moreover, we observe that the upper bound coincides with the classical diameter, $d_c(G)$, when the noise level $\xi$ goes to $0$, which corresponds to the result for DGT instances.

Note that the upper bound in Lemma \ref{lemma:xisgt new} depends only on the classical diameter of the underlying graph, $d_c(G)$, and parameters $\epsilon$, $\xi$, but not on $\tau_0(\mathcal I_G)$.
Furthermore, in contrast to Lemma \ref{lemma:xisgt diameter}, in Lemma \ref{lemma:xisgt new}, the parameter $\epsilon$ can take any value in $(0,1)$. On the flip side, where the noise parameter $\xi$ can take on any value in $(0,1)$ in Lemma \ref{lemma:xisgt diameter}, the result in Lemma \ref{lemma:xisgt new} is restricted to the case where $\xi  \in (0,1/2)$. This restriction is due largely to the limitation of our analysis. Specifically, in the proof of Lemma \ref{lemma:xisgt new}, we employ a random walk-based argument. In each step, we move one step closer to a target state with probability $1-\xi$, and one step away from the target with probability $\xi$.
Within this framework, we show that the process is able to reach a target node with high probability within a sufficiently large number of steps if the noise level $\xi<1/2$, which leads to an upper bound on $\tau_\epsilon$. When the noise level $\xi$ approaches or exceeds $1/2$, however, basic results from the theory of random walks show that the expected number of steps for a random walk to reach the target becomes infinite (see~\cite{durrett2019probability} for example). We are hopeful that improved analysis in a future work can help remedy this restriction and address the case where $\xi \geq 1/2$.

\subsection{Dynamic Energy Management with Storage}\label{subsec:dem}
In Section \ref{subsec:xisgt}, we introduced the $\xi$-SGT family of MDP instances, which is a stochastic variant of the DGT instance introduced in Section \ref{subsec:dgt}. In this subsection we provide an illustrative application that can be modeled by the $\xi$-SGT family. We consider the following model of dynamic energy management with storage.

Consider an operator and a battery with $B$ charging levels $\mathcal S = \{0,\ldots,B-1\}$, and a power parameter, $C\in\{1,\ldots,B-1\}$, representing the maximum units of charging and discharging within one time step. The state $S_t$ corresponds to the battery level at time $t$. The transition function is given by
\begin{equation}
S_{t+1} = \left(S_t+\min \left\{a_t,Y_t^S\right\}\right)_{[0,B-1]},
\end{equation}
where $(x)_{[a,b]}$ represents the projection of $x$ onto the interval $[a,b]$. Here, $Y_t^S$ is a nonnegative random variable representing the \emph{on-site renewable generation} (e.g., wind or solar) at time $t$. We will assume that $Y_t^S$ satisfies $\mathbb P\left(Y_t^S < C\right)=\beta$, and $\mathbb P\left(Y_t^S \geq C\right) = 1-\beta$. The value $\beta$ is a noise level, such that $1-\beta$ corresponds to the probability that there is enough renewable generation for the decision maker to achieve maximum per-step charge, $C$. 

The variable $a_t\in\{-C,\ldots,C\}$ represents the \emph{control} at time $t$: the decision maker may choose to sell the stored energy by setting $-C\leq a_t<0$, charge the battery by setting $0< a_t\leq C$, or hold the current battery level by setting $a_t = 0$. Any unused energy is stored in the battery, up to its capacity, $B-1$. Note that when $a_t>0$, the actual amount of energy charged to the battery is 
\begin{equation}
a_t^C = \min\left\{B-1-S_t,(a_t)^+,Y_t^S\right\}.
\end{equation}
 When $a_t<0$, the actual amount of energy sold is 
 \begin{equation}
 a_t^S = \min\{(a_t)^-,S_t\}. 
 \end{equation}
 In other words, the charging process may be impacted by the random renewable generation, while the selling and holding actions are assumed to be noiseless in this model. 

The goal of an operator is to maximize the expected total reward $V = \sum_{t=0}^{T-1} R_t\left(a_t,S_t,Y_t^R\right)$ for some reward function $R_t$, which takes values in $[r_{\min},r_{\max}]$ for some $r_{\min},r_{\max}\in\mathbb R$ with $r_{\max}-r_{\min}=\bar{r}>0$. One example of the reward function can be the operator's net revenue, defined as the difference between the revenue generated from selling energy and the charging costs. In particular, the reward function can be expanded as 
\begin{equation}
R_t \left(a_t,S_t,Y_t^R\right)= -a_t^CP_t^C+a_t^S P_t^S,
\end{equation}
where $P_t^C$, $P_t^S$ are the charging costs and the selling prices at time $t$, respectively. The prices $P_t^C$, $P_t^S$ are bounded nonnegative random variables with mean $p_t^C$, $p_t^S$, respectively. At each time step, the agent plans the control $a_t$ ahead when only the mean prices are available, but not the actual prices. 

\begin{figure}[h]
\centering
\includegraphics[scale = 0.5]{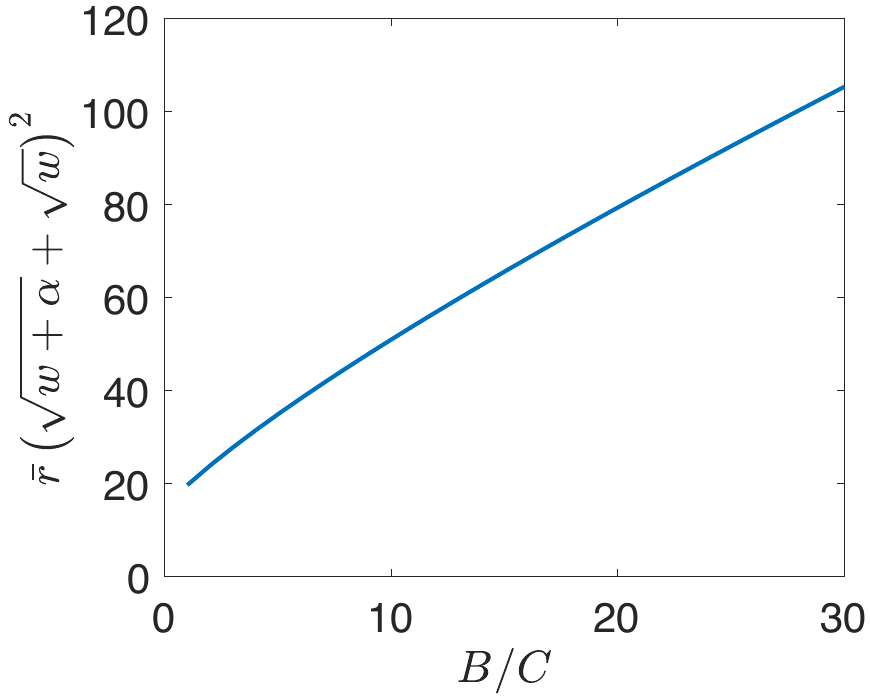}
\caption{An illustration of the upper bound in Theorem \ref{thm:energy}, with $\beta = 0.1$, $\bar{r}=1$, and $B/C$ ranging from 1 to 30.} 
\label{fig:energy_bound}
\end{figure}

It is easy to verify that the dynamic energy management system depicted above is a $\xi$-SGT problem based on $(G,\xi)$, where $\xi=\beta$ and $G$ is a connected undirected graph with $B$ vertices and a noiseless self-loop around each vertex. 
Each node in $G$ has no more than $2C+1$ edges.
In particular, for any pair of states in the state space, $s,s'\in\{0,\ldots,B-1\}$, such that $0<s'-s\leq C$, we can choose an action $a=s'-s$ such that by taking action $a$, the state transitions from $s$ to $s'$ in one step with a probability at least $1-\beta$. If $-C\leq s'-s\leq 0$, we can choose $a=s'-s$ such that by taking action $a$, the state goes from $s$ to $s'$ with probability one.  Note also that with $B$ battery levels and parameter $C$, we have
\begin{equation}\label{eq:diameter_bc}
d_c(G) \leq \frac{B}{C}+1.
\end{equation}
 With these observations in mind, by applying Theorem \ref{thm:upperbound} and Lemma \ref{lemma:xisgt new}, we have the following result that characterizes the regret of temporal concatenation in this family of problems. The proof is provided in Appendix \ref{app:thm:energy}. 
\begin{theorem}\label{thm:energy}
For a dynamic energy management system with $B$ battery levels, power parameter $C\leq B-1$, and noise level $\beta$, the regret of temporal concatenation for any initial distribution $\mu_0$ and horizon $T$ is bounded from above as follows:
\begin{equation}\label{eq:energy_bound_2}
\Delta(\mu_0,T) \leq \bar{r} \left(\sqrt{\omega+\alpha}+\sqrt{\omega}\right)^2.
\end{equation}
where we use shorthands $\alpha: = \frac{\frac{B}{C}+1}{1-2\beta}$, $\omega: = f\left(\frac{B}{C}+1,\beta\right)$, with $f(\cdot,\cdot)$ defined in Eq.~\eqref{eq:def:f}. 
\end{theorem}
Figure \ref{fig:energy_bound} contains a numerical example for the right-hand side of \eqref{eq:energy_bound_2}, illustrating its relationship with the ratio $B/C$ given fixed noise level $\beta$. To put the figure in context, let us look at  the following illustrative example.  Consider an energy-storage system with a total capacity of 36 MWh and hourly power rating of 9 MW \cite{byd2017}. Assume further that each time slot in the MDP amounts to approximately 10 minutes (e.g., the California ISO has a 5 or 15-min dispatch window for the real-time utility energy market \cite{CAISO2019}). This translates to $B=36$,  $C = 9/6  = 1.5$, and $B/C = 24$. With $\beta=0.1$ and a normalized $\bar{r}=1$, Theorem \ref{thm:energy} would suggest that the  total regret is bounded from above by  90, uniformly over all time horizons. Since there are 144 time slots in a 24-h period, this suggests that the average regret of one step incurred by temporal concatenation is at most $62.5\%$ over a one-day horizon, or $8.9\%$ for a one-week horizon.

\subsection{Simulation Results for DGT Instances}
\label{subsec:numerical}
In this section, we provide numerical examples to illustrate the trajectory of the regret of temporal concatenation, which will allow us to investigate the degree to which the theoretical results in Section \ref{sec:main_result} hold in \emph{``average'' instances} with different diameters. We also explore the performance of a generalized temporal concatenation, which temporally concatenates the policies of $K$ sub-instances, for $K\geq 2$.

We will consider the DGT instance based on a graph $G$ with finite $0$-diameter and deterministic state transition, as defined in Definition \ref{def:dgt}. Suppose the reward functions, $\{R_t\}_{t\in[T]}$, are also deterministic and depend only on the current state. In this case, each vertex in the graph $G$ corresponds to a state of the MDP instance and is associated with a reward. In graph $G$, an edge from vertex $i$ to vertex $j$ means that the system can transition from state $i$ to state $j$ within one step when an appropriate action is taken. 

We present a group of simulations where the transition function can be represented by a strongly connected directed graph with at least one self-loop.

\begin{figure}[h]
\centering
\includegraphics[scale = 1]{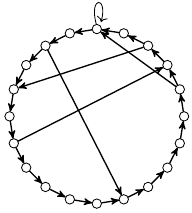}
\caption{An example of a strongly connected directed graph with at least one self-loop, $G_D^{(i)}$, with $n=20$ nodes.}
\label{fig:example}
\end{figure}

\subsubsection{Simulation set-up} We consider DGT instances depicted by directed graphs. 
{This generative model allows us to control the diameter of the instance by varying the density of the randomly added edges.}
We randomly construct $N_D = 3\times10^4$ DGT instances with $|\mathcal S| = 200$ states.
For each realization, the MDP has a deterministic transition function, which can be represented by $G_D^{(i)}\in\mathcal G_{\rm csl}$, a strongly connected directed graph including at least one self-loop, for $i\in[N_D]$. Let $p_D^{(i)} = \frac{W_D^{(i)}}{200}$, where $W_D^{(i)}\overset{\rm i.i.d.}{\sim} {\rm Unif}(0,1]$. We generate directed graphs $G_D^{(i)}$ by randomly adding edges to a ring, which is similar to the small-world model in \cite{watts1998collective}.
Specifically, the graphs $G_D^{(i)}$ are independently generated using the following steps: 

\noindent (i) construct a ring that connects all the states with edges $(1,2), (2,3)$, $\ldots,(|\mathcal S|-1,|\mathcal S|),(|\mathcal S|,1)$; 

\noindent (ii)  add a self-loop around the vertex 1; 

\noindent (iii) with probability $p_D^{(i)}$ add an edge $(j,k)$ if there is currently no edge from $j$ to $k$, for $j,k\in \mathcal S$. (If $j=k$ this will be a self-loop around $j$.) 

In the $i$th realization, the reward associated with each node $R^{(i)}(j)$ is drawn uniformly at random from the set $\{1,2,\ldots,200\}$, for $j\in\mathcal S$, which implies that the maximal reward $\bar{r}^{(i)}\leq 200$.
In Figure \ref{fig:example}, we provide an example of $G_D^{(i)}$ with 20 nodes.

Once an MDP instance is constructed, we compute the regret of temporal concatenation with a uniform initial state distribution $\mu_0 = \left(\frac{1}{|\mathcal S|},\ldots,\frac{1}{|\mathcal S|}\right)$ for different horizons $T$. Definition \ref{def:tpc} can be easily generalized to temporal concatenation with $K$ sub-instances for $K\geq 2$, which will be elaborated in the subsequent paragraph. For temporal concatenation with $K$ sub-instances with $K = 2,3,4,5$, we let $T$ vary from $K$ to 800. For each case, we run $N_D=3\times 10^4$ simulations and compute the classical diameter $d_c$ of the graph in each realization. Let $\mathcal N_d = \left\{i:\ d_c\left(G_D^{(i)}\right)=d,\ i\in[N_D]\right\}$ be the collection of all graphs generated in the simulation with classical diameter $d$. 
For realizations with the same diameter $d$, we compute the (normalized) empirical average regret of temporal concatenation for different $T$, which can be expanded as
\begin{align}\label{eq:Delta_hat}
\widehat{\Delta}(d,T)\triangleq\frac{1}{|\mathcal N_d|}\sum_{i\in\mathcal N_d}\frac{1}{\bar{r}^{(i)}}\left(\mathbb E^{\pi^*}\left[\sum_{t=0}^{T-1} R^{(i)}\left(S_t^{(i)}\right)\right]-\mathbb E^{\pi_{\rm TC}}\left[\sum_{t=0}^{T-1} R^{(i)}\left(S_t^{(i)}\right)\right]\right).
\end{align}
Here, $S_t^{(i)}$ is the state at time $t$ in the $i$th realization, $G_D^{(i)}$, while $R^{(i)}(j),\ j\in \mathcal S$ are regarded as parameters in \eqref{eq:Delta_hat}.  
Note that in $\widehat{\Delta}(d,T)$, we normalize the regret of the $i$th instance by its maximal reward $\bar{r}^{(i)}$.
For each diameter $d$, we find the (normalized) empirical maximum average regret with respect to $T$, i.e.,
\begin{equation}\label{eq:Delta_hat_max}
    \widehat{\Delta}_{\rm \max}(d) = \max_T\widehat{\Delta}(d,T),
\end{equation}
where the maximum is taken over all $T$ included in the simulation.

Now we define temporal concatenation with $K$ sub-instances $(K\geq 2)$ in an analogous way as in Definition \ref{def:tpc}. For an original instance $\mathcal I_0$, denote by $\{\mathcal I_k\}_{k\in[K]}$ the sub-instances generated by partitioning $\mathcal I_0$ into $K$ sub-instances of (approximately) equal length along the time horizon.
Let $\pi_k^*=\sol(\mathcal I_k)$, $k\in[K]$. The temporal concatenation heuristic with $K$ sub-instances generates a policy, $\pi_{\rm TC}$, by temporally concatenating optimal solutions for $\mathcal I_k$, which is analogous to \eqref{eq:tpc1}.

\begin{figure}
    \centering
    \subfloat[]{{\includegraphics[scale=.33]{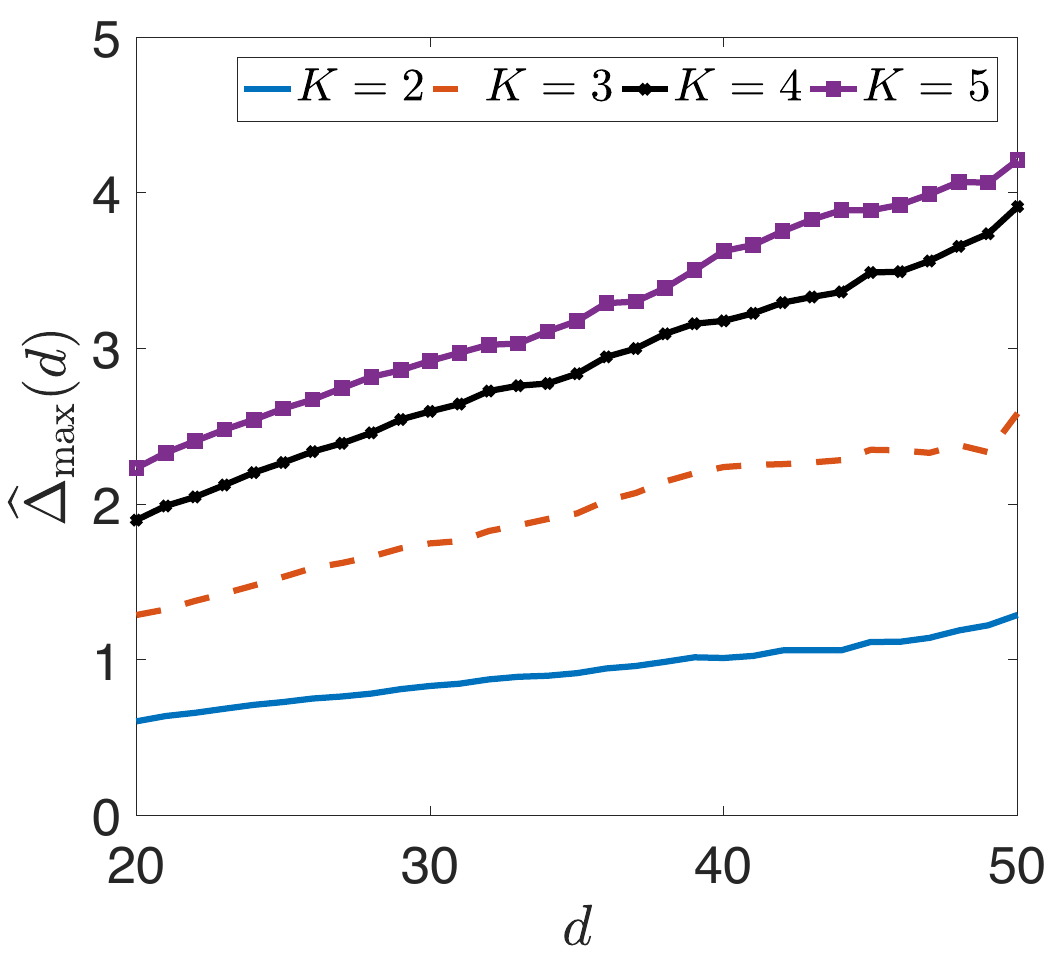}} 
    \label{fig:directed_a}
    }%
    \quad
    \subfloat[]{{\includegraphics[scale=.33]{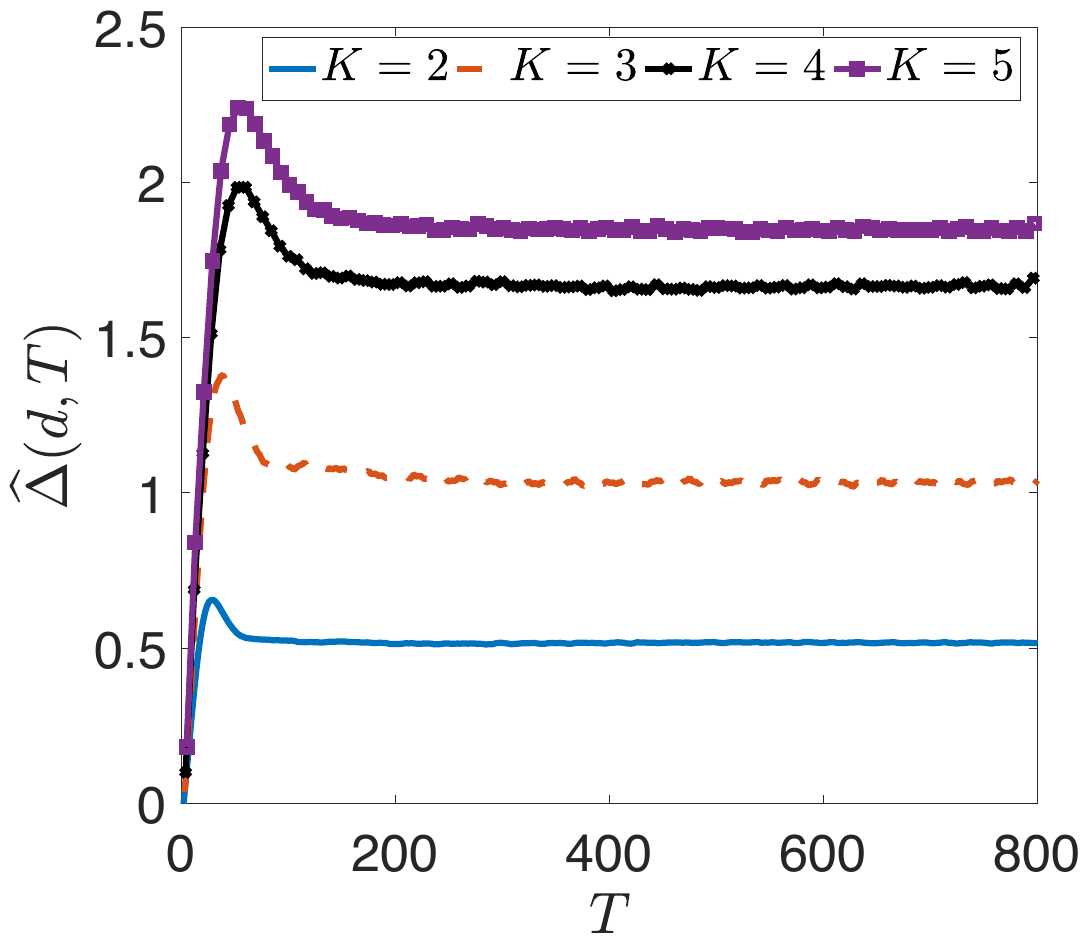} }
    \label{fig:directed_b}
    }%
    \caption{The normalized regret of DGT instances based on directed graphs for temporal concatenation with $K = 2,3,4,5$ sub-instances. (The original temporal concatenation corresponds to $K=2$.) (a)  Empirical maximum average regret $\widehat{\Delta}_{\rm max}(d)$ as a function of the diameter $d$; (b) Empirical average regret $\widehat{\Delta}(d,T)$ as a function of the horizon $T$, for a fixed diameter $d=23$. (The plots are smoothed by a 5-step moving-average filter.)}%
    \label{fig:directed}%
\end{figure}

\subsubsection{Results} Our first finding shows that the empirical maximum average regret $\widehat{\Delta}_{\rm max}(d)$ increases {linearly with respect to the diameter $d$.} As illustrated in Figure \ref{fig:directed_a}, the empirical maximum average regret $\widehat{\Delta}_{\rm max}(d)$ exhibits an increasing trend as the diameter $d$ increases from 20 to 50 for temporal concatenation with $K=2,3,4,5$ sub-instances. By Lemma \ref{lemma:connection of diameters}, for a DGT instance based on $G$, $\mathcal I_G$, the 0-diameter is bounded between $d$ and $2d$. Hence, the numerical result is consistent with Theorem \ref{thm:upperbound}, which bounds the performance regret from above by the $0$-diameter $\tau_0$ if the maximal reward is normalized to be 1. {Note that the slopes in Figure \ref{fig:directed_a}
are much smaller than $1$, which, as expected, is due to the worst-case nature of the upper bound}. From the same figure, we also see that increasing the number of sub-instances in temporal concatenation will increase the average regret. In particular, when the horizon $[T]$ is fixed, as the number of sub-instances, $K$, increases, the length of each sub-instance decreases. Shorter sub-instances will more likely lead to overly short-sighted policies, which impede the performance of temporal concatenation. 

The second finding suggests that for a fixed diameter $d$, as $T$ grows, the empirical average regret $\widehat{\Delta}(d,T)$ first increases, then decreases after reaching a peak, and finally stabilizes when $T$ is sufficiently large. This trend is illustrated in Figure \ref{fig:directed_b}. Intuitively, when $T$ starts growing from zero, temporal concatenation starts to incur performance regret. Since the temporal concatenation policy is sub-optimal, the regret becomes larger with more time steps. When $T$ is sufficiently large, however, the regret no longer increases. An intuitive explanation is that the temporal concatenation policy and the optimal policy become similar when the length of a sub-instance is sufficiently large, which causes the regret to start decreasing in this region. It remains an interesting open problem {for} finding the minimum horizon $T$ beyond which the average regret starts to decrease.

\subsection{Run-time Reduction by Temporal Concatenation}\label{subsec:runtime}
In this subsection, we conduct additional numerical simulations to assess the benefit of run-time reduction from using temporal concatenation. We apply the classic value iteration algorithm to solve DGT instances, and compare the run-time of the following two cases: 

\begin{enumerate}
    \item sequentially computing the optimal policy of the original instance; 
    \item using temporal concatenation, where we employ the built-in Matlab command \texttt{spmd} to solve the two sub-instances in parallel on a multi-core processor. 
\end{enumerate}

We  run all experiments on a standard multi-core desktop computer; the specifications of the environment are described in detail in Appendix \ref{app:sys_spec}. We construct DGT models with the same setting as described in Section \ref{subsec:numerical}, except that here we vary the number of states and time horizon length. We consider the cases that the number of states $|\mathcal S|=1000, 2000, 3000$, and the time horizon $T = 50, 500, 12000$, respectively. For each pair of $|\mathcal S|$ and $T$, we randomly generate 30 instances and present their respective run-time in Tables \ref{tab:n1000}, \ref{tab:n2000}, and \ref{tab:n3000}.
We will denote by $T_{\rm seq}$  the run-time of solving the original problem sequentially.  $T_{\rm tc}$ denotes the time taken by the temporal concatenation method when run with Matlab's native parallel computation framework. The ratio between the two values $\eta = T_{\rm tc}/T_{\rm seq}$ is therefore a metric of interest, showing the multiplicative speed-up obtained by temporal concatenation. For a value $x$, we will use $\overline{x}$ to denote its empirical mean,  and  $\sigma_x$ its standard deviation.

We can see from the results that the temporal concatenation heuristic reduces the run-time thanks to parallelism. We also note that the time reduction is not exactly 50\% as one might expect, and this is largely due to the fact that there is a non-trivial overhead for Matlab to initiate a parallel computation instance. That being said, the reduction is still fairly significant across the board, and gets closer to 50\% when the run-time of the original problem is sufficiently long (e.g., for the cases of $|\mathcal S|=3000$, $T=500, 12000$).  

Interestingly, we notice that the value of $|\mathcal S|$ seems more significant in determining the ratio $\eta$ than the horizon $T$. Suppose we fix $|\mathcal S|$. In this case, increasing $T$ from 50 to 500 reduces $\eta$ notably, while further increasing its value from 500 to 12000 seems not very influential. On the other hand, if we fix $T$ and increase $|\mathcal S|$ from 1000 to 2000 and then to 3000, the value of $\eta$ decreases remarkably. We conjecture that this is due to the fact that the run-time of each sub-instance scales quadratically in $|\mathcal S|$ but only linearly in $T$. Therefore, as the two sub-instances become larger, its computation time would tend to dwarf Matlab's computational overhead in setting up the parallel computation instance, leading to a more favorable ratio of speed-up  $\eta$. 

Admittedly, what we have here is a relatively simple proof-of-concept with off-the-shelf software, and we expect that a more optimized implementation of the parallel computation, possibly on distinct physical machines, would further decrease the run-time.

\begin{table}[ht]
\centering
\subfloat[][$|\mathcal S| = 1000$]{
\begin{tabular}{|l@{\hskip\arraycolsep}l@{\hskip\arraycolsep}l@{\hskip\arraycolsep}l|}
\hline
$T$ & $\overline{T_{\rm seq}}$ ($\sigma_{T_{\rm seq}}$)& $\overline{T_{\rm tc}}$ ($\sigma_{T_{\rm tc}}$) & $\overline{\eta}$ ($\sigma_\eta$) \\
\hline
50 & 0.292 (0.007) & 0.280 (0.020) & 0.958 (0.060)\\
500 & 2.881 (0.058) & 2.065 (0.117) & 0.717  (0.042) \\
12000 & 68.494 (1.542) & 48.482 (2.182)  & 0.708 (0.035)\\
  \hline
\end{tabular}
\label{tab:n1000}
}
\\
\subfloat[][$|\mathcal S|=2000$]{
\begin{tabular}{|l@{\hskip\arraycolsep}l@{\hskip\arraycolsep}l@{\hskip\arraycolsep}l|}
\hline
$T$ & $\overline{T_{\rm seq}}$ ($\sigma_{T_{\rm seq}}$) & $\overline{T_{\rm tc}}$ ($\sigma_{T_{\rm tc}}$) & $\overline{\eta}$ ($\sigma_\eta$) \\
\hline
50 & 1.463 (0.033) & 0.989 (0.034) & 0.676 (0.022) \\
500 & 14.747(0.149) & 8.883 (0.248) & 0.602 (0.018) \\
12000 & 349.295 (5.910) & 212.220 (10.334) & 0.608 (0.020)\\
  \hline 
\end{tabular}
\label{tab:n2000}
}
\\
\subfloat[][$|\mathcal S|=3000$]{
\begin{tabular}{|l@{\hskip\arraycolsep}l@{\hskip\arraycolsep}l@{\hskip\arraycolsep}l|}
\hline
$T$ & $\overline{T_{\rm seq}}$ ($\sigma_{T_{\rm seq}}$) & $\overline{T_{\rm tc}}$ ($\sigma_{T_{\rm tc}}$)& $\overline{\eta}$ ($\sigma_\eta$) \\
\hline
50 & 3.514 (0.026) & 2.135 (0.061) & 0.608 (0.017) \\
500 & 35.381 (0.143) & 19.633 (0.327) & 0.555 (0.009) \\
12000 & 837.234 (15.070) & 472.133 (14.640) & 0.564 (0.015)\\
  \hline
\end{tabular}
\label{tab:n3000}
}
\caption{Computation time comparison between sequential value iteration and temporal concatenation.}
\label{tab:runtime}
\end{table}

\subsection{Simulation Results for the GARNET MDP Model}\label{subsec:garnet}

In this subsection, we provide numerical simulations for a more widely studied family of MDP instances introduced in \cite{bhatnagar2009natural}, known as the GARNET model. 
In this generative model of MDP, a branching factor $B$ determines the transition kernels. When sampling an MDP instance, under each action, every state is randomly assigned $B$ possible next states, while the probability of transitioning to each of them is also randomly drawn. Here, we construct GARNET MDP models with $|\mathcal S|=200$ states, $|\mathcal A|=3$ actions, and let $B$ range from 3 to 15. For each choice of $B$, we construct $N_G=5000$ examples independently and record their average regrets. To make sure that the MDP is aperiodic, we introduce a self-loop around each state for every action, i.e., each state has one edge to itself and $B-1$ edges to other states under each action. The reward functions $R_t$ are deterministic, which are determined by the current state and the chosen action. For each pair of state and action, the corresponding reward is drawn uniformly at random from 1 to 200. In this simulation, we focus on the regret of temporal concatenation with two sub-instances.

\begin{figure}[h!]
    \centering
\includegraphics[width=0.6\textwidth]{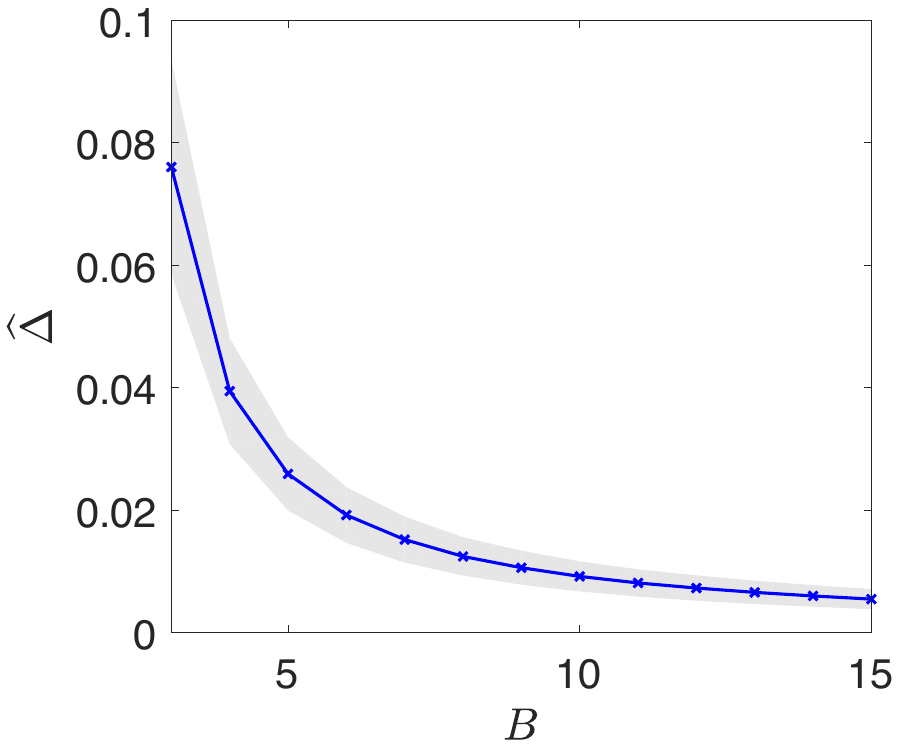}
    \caption{The normalized regret of temporal concatenation with two sub-instances for the GARNET MDP model with $|\mathcal S|$ =200, $T=200$, $|\mathcal A|=3$, $B=3,4,\ldots,15$. }
    \label{fig:garnet}
\end{figure}

In Figure \ref{fig:garnet}, we present the simulation results. 
Here, $\widehat{\Delta}$ is the empirical average regret normalized by the maximal reward, defined as follows:
\begin{equation}
    \widehat{\Delta} = \frac{1}{N_G}\sum_{i=1}^{N_G} \frac{1}{\bar{r}^{(i)}} \left(\mathbb E^{\pi^*}\left[\sum_{t=0}^{T-1} R_t^{(i)}\left(a_t^{(i)},S_t^{(i)}\right)\right]-\mathbb E^{\pi_{\rm TC}}\left[\sum_{t=0}^{T-1} R_t^{(i)}\left(a_t^{(i)},S_t^{(i)}\right)\right]\right),
\end{equation}
where $R_t^{(i)}$, $a_t^{(i)}$, $S_t^{(i)}$, $\bar{r}^{(i)}$ are, respectively, the reward functions, actions, states, and maximal reward in the $i^{\rm th}$ realization of the GARNET model, $i\in [N_G]$. 

 The shaded area in the figure depicts the area within one empirical standard deviation of the mean. Overall, the performance is favorable.  We see that the regret in the GARGET model is substantially smaller than that of DGT (Figure \ref{fig:directed}), and it decreases even more as $B$ becomes large. We suspect that the random and uniform nature with which GARGET generates transition kernels contributes to the resulting MDP having a relatively small $\epsilon$-diameter, even when $B$ is moderate, and the diameter becomes even smaller as the number of neighboring states, $B$, increases.

\section{Conclusion}\label{sec:conclusion}
In this paper, we propose and analyze a heuristic architecture, temporal concatenation, for speeding up existing MDP algorithms when solving a finite-horizon Markov decision process. Temporal concatenation decomposes the problem over the time horizon into smaller sub-problems and subsequently concatenating their optimal policies. 

Using a notion of $\epsilon$-diameter, we provide upper bounds that show, when the underlying MDP instance admits a bounded $\epsilon$-diameter the regret of temporal concatenation is bounded and \emph{independent} of the length of the horizon. Conversely, we provide lower bounds by showing that, for any finite diameter, there exist MDP instances for which the regret upper bound is tight for all sufficiently large horizons.

At the high level, we aim to explore an alternative approach for solving large-scale MDPs: instead of creating new algorithms from scratch, we may be able to leverage existing MDP algorithms in creative ways to harness additional performance gains. The present paper takes a first step towards this direction by decomposing the problem  along the time axis. There is a number of interesting directions for future work. While we have demonstrated that the $\epsilon$-diameter of an MDP has a substantial impact on the performance of temporal concatenation, it remains a challenge in general to compute such diameter for a given MDP. Understanding how to obtain sharp bounds for, or numerically compute, the $\epsilon$-diameter can be an interesting direction of future research.
The theory of our present paper has mostly focused on dividing the MDP into two equal-sized sub-instances, while in general, one may consider $K\geq 2$ sub-instances, or study the case where the sizes of the sub-instances may even vary. It would be interesting to understand how best to choose the number and the sizes of the sub-instances, especially when the original MDP is time-inhomogeneous.

In addition, it will be interesting to explore connections between temporal concatenation and other approximation heuristics aimed at reducing the complexity of solving an MDP. For instance, if we further assume that the MDP is time-homogeneous (transitions and rewards do not vary with time) and that the time horizon is very large, then a natural alternative would be to directly deploy a stationary optimal policy associated with the average-reward version of the MDP, which can be calculated efficiently using a horizon-independent linear program. The regret of such an approach would likely depend on whether the steady-state optimal policy could quickly reach its stationary distribution from an arbitrary initial condition, and it would be interesting to understand, for instance, whether the latter is related to the notion of diameter  studied in this paper.

Finally, in a broader sense, the present work only explores decomposing an MDP along the time axis, and it would be interesting to explore the efficiency of other forms of decomposition architectures, such as those that operate through states. 
\section*{Acknowledgements}
We thank the anonymous referees for their comments and feedback.

\bibliographystyle{acm}
\bibliography{tc_new}

\newpage
\appendix
\section*{Appendices}

\section{Proofs of Main Results}\label{sec:proof}

This section is devoted to the proofs of the main results. Before delving into the details, we first provide a high-level overview of the key ideas.

\emph{Upper bounds} (Section \ref{sec:proof_upper}). For Theorem \ref{thm:upperbound}, observe that the temporal concatenation heuristic by construction achieves optimal total expected reward during the first sub-instance, $\calI_1$. The problem arises, however, if acting greedily during  $\calI_1$ would result in the system being in a disadvantageous state at the beginning of the second sub-instance, $\calI_2$, thus leading to a large regret. Our analysis for the regret upper bound in Theorem \ref{thm:upperbound} will therefore focus on the dynamics of temporal concatenation during $\calI_2$. To this end, we will employ a coupling argument, by bounding temporal concatenation's regret from above using that of a carefully constructed, and likely strictly sub-optimal, ``fictitious'' policy, $\tilde{\pi}$, during $\calI_2$. The policy $\tilde \pi$ consists of multiple phases of length approximately $\tau_\epsilon$. In the $k$th phase, it aims to reduce the (total variation) distance with the overall optimal policy $\pi^*$ over the course of $\tau_\epsilon$ steps. Using an argument based on recursion, we show that in the $k$th phase this policy incurs a regret that is up to $\epsilon^{k-1}\bar{r}\tau_\epsilon$. This will in turn allow us to show that the regret of $\tilde{\pi}$ incurred during the second phase is small.

\emph{Lower bounds} (Sections \ref{sec:proof_lower_fin_Diam}). For the lower bound in Theorem \ref{thm:lowerb_finite_diam}, we build on the insights gathered from the proof of Theorem \ref{thm:upperbound} to generate worst-case MDP instances. The main idea is to construct instances in such a way that during the first sub-instance, the temporal concatenation heuristic is guaranteed to be lured by some \emph{small} short-term rewards and end up in a ``bad'' subset of the state space, from which it will suffer large losses in the second half of the time horizon compared to the optimal policy.

\subsection{Proof of Theorem \ref{thm:upperbound}}
\label{sec:proof_upper}

Recall that $\mu^\pi_t$ is the distribution of the state $S_t$ induced by a policy $\pi$.  For simplicity of notation, we will write $\mu_t$ in place of $\mu^\pi_t$ when there is no ambiguity, and use the shorthand: 
\begin{equation}
\mu_t^{\pi^*}\triangleq \mu_t^*,  \quad \mbox{ and }  \mu_t^{\tc}\triangleq \mu_t^{\rm TC}. 
\end{equation}
 We also define the cumulative rewards from time $t_1$ to $t_2$ as 
\begin{equation}
\tilde{V}(t_1,t_2) \triangleq \sum_{t=t_1}^{t_2} R_t(a_t,S_t,Y^R_t), \quad t_1 < t_2.
\end{equation}

The next lemma is the key technical result.  
Recall from \eqref{eq:defReg} that, for a given instance, the total expected reward of a policy depends on the initial distribution. In Lemma \ref{lem:upper_bound}, we provide an upper bound on the performance difference of the optimal policy when the system starts from two different initial distributions.
\begin{lemma}
\label{lem:upper_bound}
Fix an instance $\calI$ with horizon $[T]$. Fix distributions $\mu_0,\nu_0\in\mathcal P$. Let $\pi^*$ be the optimal policy for the instance $\mathcal I$. If there exists $\epsilon>0$ such that $\tau_\epsilon(\calI)\leq T$, the difference in total expected reward under $\pi^*$ between the cases where the initial distribution is $\mu_0$ versus $\nu_0$ is bounded from above as follows:
\begin{equation}
       \left|V(\calI, \pi^*, \mu_0)-V(\calI, \pi^*, \nu_0)\right|\leq \frac{\bar{r}\tau_\epsilon(\mathcal I)}{1-\epsilon}.
        \end{equation}
\end{lemma}
We first prove Lemma \ref{lem:upper_bound}, using a coupling argument.
For state $s\in\mathcal S$, starting time $t\in[T]$, and policy $\pi$, we define the value function as follows:
\begin{equation}
V_t^\pi(s) =  \E^{\pi}\left[\tilde{V}(t,T-1)\Big|S_t=s\right].
\end{equation}
 For the instance $\mathcal I$ and policy $\pi^*$, the total expected reward for initial distribution $\mu$ is 
\begin{equation}
V(\mathcal I, \pi^*,\mu) = \sum_{s\in\mathcal S} \E^{\pi^*}\left[\tilde{V}(0,T-1)\Big|S_0=s\right]\mu(s) = \sum_{s\in\mathcal S}V_0^{\pi^*}(s)\mu(s).
\end{equation}
Fixing the policy $\pi^*$, the difference in total expected rewards under $\pi^*$ but starting with two initial distributions, $\mu_0$ and $\nu_0$, can be expanded as:
\begin{align}
\left|V(\mathcal I, \pi^*,\mu_0)-V(\mathcal I,\pi^*,\nu_0)\right|= \left|\sum_{s\in\mathcal S} V_0^{\pi^*}(s)(\mu_0(s)-\nu_0(s))\right|.\label{pf:6}
\end{align}
Without loss of generality, suppose that  
\begin{equation}\label{eq:wlog}
V(\mathcal I, \pi^*,\mu_0)\geq V(\mathcal I,\pi^*,\nu_0).
\end{equation}
Now we provide an upper bound on the difference in total expected reward by introducing a ``fictitious'' policy $\tilde{\pi}$. Suppose $\tau_\epsilon(\mathcal I)\leq T$ for some $\epsilon>0$.
Recall that by the definition of $\epsilon$-diameter, there exists a policy $\tilde{\pi}_{\rm ap}$ such that starting from $\nu_0$, the state distribution at time $\tau_\epsilon(\mathcal I)$ under $\tilde{\pi}_{\rm ap}$, which is denoted by $\tilde{\nu}_{\tau_\epsilon(\mathcal I)}$, satisfies 
\begin{equation}
\delta_{\rm TV} (\mu^*_{\tau_\epsilon(\mathcal I)},\tilde{\nu}_{\tau_\epsilon(\mathcal I)})\leq \epsilon,
\end{equation}
where $\mu^*_{\tau_\epsilon(\mathcal I)}$ is the state distribution at time $\tau_\epsilon(\mathcal I)$ starting from $\mu_0$ under policy $\pi^*$.
The policy $\tilde{\pi}$ is defined as follows: for time $t\in 0 \to \tau_\epsilon(\mathcal I)-1$, let $\tilde{\pi} = \tilde{\pi}_{\rm ap}$; for time $t\in \tau_\epsilon(\mathcal I)\to T-1$, let $\tilde{\pi} = \pi^*$.
Note that $\tilde{\pi}$ is sub-optimal compared to $\pi^*$, and we have
\begin{equation}\label{eq:tilde_pi}
V(\mathcal I,\pi^*,\nu_0)\geq V(\mathcal I,\tilde{\pi},\nu_0).
\end{equation}
Recall that the reward function $R_t$ takes values in $[0,\bar{r}]$. We have that
\begin{align}
\left|V(\mathcal I, \pi^*,\mu_0)-V(\mathcal I,\pi^*,\nu_0)\right| 
= &V(\mathcal I, \pi^*,\mu_0)-V(\mathcal I,\pi^*,\nu_0)\label{pf:10-1} \\
\leq & V(\mathcal I, \pi^*,\mu_0)-V(\mathcal I,\tilde{\pi},\nu_0)\label{pf:10-2} \\
\leq & \bar{r}\tau_\epsilon(\mathcal I)+ \sum_{s\in\mathcal S} V_{\tau_\epsilon(\mathcal I)}^{\pi^*}(s)(\mu^*_{\tau_\epsilon(\mathcal I)}(s)-\tilde{\nu}_{\tau_\epsilon(\mathcal I)}(s)),
\label{pf:10-3}
\end{align}
where \eqref{pf:10-1} follows from \eqref{eq:wlog}, \eqref{pf:10-2} from \eqref{eq:tilde_pi}, and \eqref{pf:10-3} from $R_t\leq \bar{r}$.

For $s\in\mathcal S$, let $\omega(s) = \min \{\mu^*_{\tau_\epsilon(\mathcal I)}(s),\tilde{\nu}_{\tau_\epsilon(\mathcal I)}(s)\}$. 
Define $\epsilon_0 = \sum_{s\in\mathcal S}\left(\tilde{\nu}_{\tau_\epsilon(\mathcal I)}(s)-\omega(s)\right)$. Note that $ \sum_{s\in\mathcal S}\left({\mu}^*_{\tau_\epsilon(\mathcal I)}(s)-\omega(s)\right) = \sum_{s\in\mathcal S}\left(\tilde{\nu}_{\tau_\epsilon(\mathcal I)}(s)-\omega(s)\right)=\epsilon_0$.

Let $\mu^-_{\tau_\epsilon(\mathcal I)}(s) = \frac{\mu^*_{\tau_\epsilon(\mathcal I)}(s)-\omega(s)}{\epsilon_0}$, and $\nu^-_{\tau_\epsilon(\mathcal I)}(s) = \frac{\tilde{\nu}_{\tau_\epsilon(\mathcal I)}(s)-\omega(s)}{\epsilon_0}$. By the definition of total variation, we have $\epsilon_0\leq \epsilon$.
Note that $\mu^-_{\tau_\epsilon(\mathcal I)},\nu^-_{\tau_\epsilon(\mathcal I)}\geq 0$, and $\sum_{s\in\mathcal S}\mu^-_{\tau_\epsilon(\mathcal I)}(s)=1$,  $\sum_{s\in\mathcal S}\nu^-_{\tau_\epsilon(\mathcal I)}(s)=1$. Hence, $\mu^-_{\tau_\epsilon(\mathcal I)}$, $\nu^-_{\tau_\epsilon(\mathcal I)}$ are probability distributions. 

Then 
\begin{align}
\left|V(\mathcal I, \pi^*,\mu_0)-V(\mathcal I,\pi^*,\nu_0)\right|  
= & \left|\sum_{s\in\mathcal S} V_0^{\pi^*}(s)(\mu_0(s)-\nu_0(s))\right|\nonumber\\
\leq & \bar{r}\tau_\epsilon(\mathcal I)+ \sum_{s\in\mathcal S} V_{\tau_\epsilon(\mathcal I)}^{\pi^*}(s)(\mu^*_{\tau_\epsilon(\mathcal I)}(s)-\tilde{\nu}_{\tau_\epsilon(\mathcal I)}(s))\nonumber\\
= & \bar{r}\tau_\epsilon(\mathcal I)+\epsilon_0 \sum_{s\in\mathcal S} V_{\tau_\epsilon(\mathcal I)}^{\pi^*}(s)(\mu^-_{\tau_\epsilon(\mathcal I)}(s)-\nu^-_{\tau_\epsilon(\mathcal I)}(s))\label{pf:11-1}\\
\leq &   \bar{r}\tau_\epsilon(\mathcal I)+\epsilon_0 \left|\sum_{s\in\mathcal S} V_{\tau_\epsilon(\mathcal I)}^{\pi^*}(s)(\mu^-_{\tau_\epsilon(\mathcal I)}(s)-\nu^-_{\tau_\epsilon(\mathcal I)}(s))\right|\nonumber\\
\leq &  \bar{r}\tau_\epsilon(\mathcal I)+\epsilon \left|\sum_{s\in\mathcal S} V_{\tau_\epsilon(\mathcal I)}^{\pi^*}(s)(\mu^-_{\tau_\epsilon(\mathcal I)}(s)-\nu^-_{\tau_\epsilon(\mathcal I)}(s))\right|.\label{pf:11-2}
\end{align}
Here \eqref{pf:11-1} follows from the definition of $\mu^-_{\tau_\epsilon(\mathcal I)}$ and $\nu^-_{\tau_\epsilon(\mathcal I)}$, and \eqref{pf:11-2} from $\epsilon_0\leq \epsilon$.
Let $N = \left\lfloor \frac{T}{\tau_\epsilon(\mathcal I)}\right\rfloor$. For $k = 0,1,\ldots,N-1$, starting from time $t = k\tau_\epsilon(\mathcal I)$, we can use the same argument to derive the following inequality:
\begin{equation}\label{eq:recursion}
\left|\sum_{s\in\mathcal S} V_{k\tau_\epsilon(\mathcal I)}^{\pi^*}(s)(\mu^-_{k\tau_\epsilon(\mathcal I)}(s)-\nu^-_{k\tau_\epsilon(\mathcal I)}(s))\right| \leq  \bar{r}\tau_\epsilon(\mathcal I)+\epsilon \left|\sum_{s\in\mathcal S} V_{(k+1)\tau_\epsilon(\mathcal I)}^{\pi^*}(s)(\mu^-_{(k+1)\tau_\epsilon(\mathcal I)}(s)-\nu^-_{(k+1)\tau_\epsilon(\mathcal I)}(s))\right|,
\end{equation}
where $\mu^-_{k\tau_\epsilon(\mathcal I)}(s)$ and $\nu^-_{k\tau_\epsilon(\mathcal I)}(s)$ are defined in the same way as $\mu^-_{\tau_\epsilon(\mathcal I)}(s)$ and $\nu^-_{\tau_\epsilon(\mathcal I)}(s)$.
Note also that 
\begin{equation}\label{eq:tail_condition}
V_{N\tau_\epsilon(\mathcal I)}^{\pi^*}(s)\leq \bar{r}\tau_\epsilon(\mathcal I).
\end{equation}
 With \eqref{eq:recursion} and \eqref{eq:tail_condition}, we have 
\begin{equation}
\left|V(\mathcal I, \pi^*,\mu_0)-V(\mathcal I,\pi^*,\nu_0)\right|  \leq \bar{r}\tau_\epsilon(\mathcal I)\left(1+\epsilon+\ldots+\epsilon^N\right)\leq \frac{\bar{r}\tau_\epsilon(\mathcal I)}{1-\epsilon}.
\end{equation}
This completes the proof of Lemma \ref{lem:upper_bound}.

Lemma \ref{lem:upper_bound} suggests that starting from two different initial distributions, $\mu_0$ and $\nu_0$, the optimal policy $\pi^*$ is guaranteed to have similar performances when the $\epsilon$-diameter is small.  We now prove Theorem \ref{thm:upperbound} using Lemma \ref{lem:upper_bound}. First, for any initial distribution $\mu_0$, we can expand the regret of temporal concatenation as the sum of the regrets incurred during the first and second sub-instance, separately: 
\begin{align}
  \Delta(\mathcal I_0,\mu_0) =& V(\calI_0, \pi^*, \mu_0) - V(\calI_0, \tc, \mu_0) \nln
= & \left(V(\calI_1, \pi^*, \mu_0) - V(\calI_1, \pi^*_1, \mu_0)\right)  \nln
&+ \left( V\left(\calI_2, \pi^*, \mu^*_{\frac{T}{2}}\right) - V \left(\calI_2, \pi^*_2, \mu^{\rm{TC}}_{\frac{T}{2}} \right) \right), 
\label{eq:delta1}
\end{align}
where, with a slight abuse of notation, we use $ V\left(\calI_1, \pi^*, \mu_0\right)$ to denote the total expected reward from applying the policy $\pi^*$ during the first sub-instance. 
Note that during the first $T/2$ steps, the original optimal policy, $\pi^*$, does not necessarily maximize the reward for this sub-instance, because it aims at maximizing the overall reward of $\calI_0$. Hence, for this sub-instance only, the temporal concatenation method is performing better than, or equally to, the original optimal policy, i.e., the first term in \eqref{eq:delta1} satisfies: 
\begin{equation}
    V(\calI_1, \pi^*, \mu_0) - V(\calI_1, \pi^*_1, \mu_0)\leq 0.
    \label{eq:regret_inst_1}
\end{equation}

We now bound the second term in \eqref{eq:delta1}. Suppose for some $\epsilon>0$, we have $\tau_\epsilon(\calI_0)\leq T/2$. 
Note that both $\pi^*$ and $\pi^*_2$ achieve the optimal performance for the second sub-instance $\mathcal I_2$, we have 
\begin{equation}
V \left(\calI_2, \pi^*_2, \mu^{\rm{TC}}_{\frac{T}{2}} \right) = V \left(\calI_2, \pi^*, \mu^{\rm{TC}}_{\frac{T}{2}} \right).
\end{equation}
Hence, using Lemma \ref{lem:upper_bound}, we have that 
\begin{align}
V\left(\calI_2, \pi^*, \mu^*_{\frac{T}{2}}\right) - V \left(\calI_2, \pi^*_2, \mu^{\rm{TC}}_{\frac{T}{2}} \right) 
= &  V\left(\calI_2, \pi^*, \mu^*_{\frac{T}{2}}\right) - V \left(\calI_2, \pi^*, \mu^{\rm{TC}}_{\frac{T}{2}} \right)\nonumber\\
\leq & \frac{\bar{r} \tau_\epsilon(\mathcal I_2)}{1-\epsilon}\label{pf:12-1}\\
\leq & \frac{\bar{r} \tau_\epsilon(\mathcal I_0)}{1-\epsilon}\label{pf:12-2},
\end{align}
where \eqref{pf:12-1} is derived by applying Lemma \ref{lem:upper_bound} to the second sub-instance $\mathcal I_2$ for initial distributions $\mu_{\frac{T}{2}}^*$ and $\mu_{\frac{T}{2}}^{\rm TC}$, and \eqref{pf:12-2} follows from the fact that $\mathcal I_2$ is a sub-instance of $\mathcal I_0$, which leads to $\tau_\epsilon(\mathcal I_2)\leq \tau_\epsilon(\mathcal I_0)$.

To complete the proof, we substitute the regret upper bounds for the first \eqref{eq:regret_inst_1} and second \eqref{pf:12-2} sub-instances into \eqref{eq:delta1}, and obtain
\begin{align}
V(\mathcal I_0,\pi^*,\mu_0)-V(\mathcal I_0,\tc,\mu_0)    \leq \frac{\bar{r} \tau_\epsilon(\mathcal I_0)}{1-\epsilon}.
    \label{eq:tildepidiff}
\end{align}
This completes the proof of Theorem \ref{thm:upperbound}. \qed

\subsection{Proof of Theorem \ref{thm:lowerb_finite_diam}} 
\label{sec:proof_lower_fin_Diam}
We now prove Theorem \ref{thm:lowerb_finite_diam} by constructing a family of MDP instances and showing that temporal concatenation suffers the regret given in the theorem on problems from this family. The key intuition is that the instances can be constructed in such a way that the temporal concatenation heuristic will be led astray by some small short-term rewards in the first half of the horizon and end up in a bad subset of the state space, from which it will suffer large losses in the second half of the time horizon compared to the optimal policy. 

Fix $d_0\geq 5$, and let $k=d_0-2$. Consider the MDP instance depicted in Figure \ref{fig:prop_large_regret}. The state space has $\vert \mathcal S\vert=k+2$ elements and we have $d_0 = k+2$. The transition function is deterministic. In states $d_1,\ldots, d_k$, and $e$, the agent can choose between two actions, such that the system either stays in the same state or goes to the next state to the right. In state $f$, the system will always go to state $d_1$ in the next step. For state $s \in\{d_1,\ldots,d_k\}$, the reward function $R_t(a,s,y)=0$, for all $t$, $a\in\calA$, and $y \in\calY^R$.  For state $s=e$ and $f$, the reward $R_t(a,s,y)$ is always equal to $\bar{r}-\sigma_1$, and $\bar{r}$, respectively, where $\sigma_1\in(0,\sigma)$ is a constant to be specified subsequently. 

We first verify that the setting has a finite diameter $\tau_0(\calI_0)=d_0=k+2$. 
Suppose the system starts from an initial distribution $\nu\in\mathcal P$ and we try to reach another distribution $\nu'\in\mathcal P$. 
Consider the following policy:

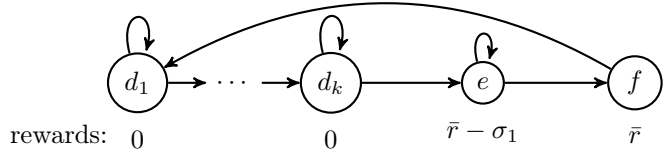
\begin{figure}
\centering
\begin{tikzpicture}[->,>=stealth',auto,node distance=2cm,
  thick,main node/.style={circle,draw,font=\sffamily}]
  \node[draw=none] (0) {};
  \node[below = 0.3cm of 0]{rewards:};
  \node[main node] (1) [right = 0.5cm of 0]{$d_1$};
  \node[below = 0.1cm of 1]{$0$};
  \node[draw=none,  right=0.5cm of 1]   (1-2) {$\cdots$};
  \node[main node] (2) [right=0.5cm of 1-2] {$d_k$};
  \node[below = 0.1cm of 2]{$0$};
  \node[main node] (3) [right of=2] {$e$};
  \node[below = 0.1cm of 3]{$\bar{r}-\sigma_1$};
  \node[main node] (4) [right of=3] {$f$};
  \node[below = 0.1cm of 4]{$\bar{r}$};
  \path[every node/.style={font=\sffamily\small}]
    (1) edge node [right] {} (1-2)
    (1-2) edge node {} (2)
    (2) edge node [right] {} (3)
    (3) edge node [right] {} (4)
    (4) edge[bend right] node [left] {} (1)
        (1) edge[loop above] node {} (1)
        (2) edge[loop above] node {} (2)
    (3) edge[loop above] node {} (3);
\end{tikzpicture}
\caption{An MDP instance with bounded diameter and large performance regret.}
\label{fig:prop_large_regret}
\end{figure}

Stage 1: If the initial state is in $\{d_1,\ldots,d_k, e\}$, stay for one step; if the initial state is $f$, go to $d_1$ in the first step. Hence, Stage 1 takes 1 step.   

Stage 2: Starting from one of the states in $\{d_1,\ldots,d_k, e\}$, the agent reaches the state distribution $\nu'$ after another $k+1$ steps. Note that starting from any state in $\{d_1,\ldots,d_k, e\}$, the system can reach any state $s\in\mathcal S$ using $k+1$ steps by first staying at the current state for an appropriate number of steps and then moving forward to reach the target state. We refer to this stay-and-move process as a $(k+1)$-path to state $s$. The agent can thus employ the following randomized policy: starting at state $s_0\in\{d_1,\ldots,d_k,e\}$, with probability $\nu'(s)$ the agent chooses to take the $(k+1)$-path to state $s$, for $s\in\mathcal S$. Stage 2 takes $k+1$ steps. 

Using the policy described above, we can reach any distribution $\nu'$ at time $t = k+2$ starting from any initial distribution $\nu$. Hence, we have shown that the diameter satisfies
\begin{equation}
\label{pf:3.3.1}
    \tau_0(\calI_0)\leq k+2. 
\end{equation}

We now establish a lower bound for $\tau_0(\calI_0)$. Suppose the initial distribution is concentrated on state $f$, i.e., $\nu(f)=1$. At time $t=0,1$, the state will be deterministically $f$ and $d_1$, respectively. Hence, in order to reach a distribution $\nu'$ with $\nu'(d_1)=\nu'(f)=0.5$, it takes at least another $k+1$ steps. Then, 
\begin{equation}\label{pf:3.3.2}
\tau_0(\calI_0)\geq k+2.
\end{equation}
In light of \eqref{pf:3.3.1} and \eqref{pf:3.3.2}, we conclude that 
$
\tau_0(\calI_0) = k+2 = d_0. 
$

We now consider the regret. Suppose the initial state is deterministically $d_1$. Recall that $T> 2d_0+2 = 2k+6$. For the optimal policy, the agent will first go to state $e$, stay at $e$ until time $T-2$, and finally go to $f$ at time $T-1$. Hence, the total reward of the optimal policy is 
\begin{equation}
V(\mathcal I_0,\pi^*,\mu_0)=(T-k)\left(\bar{r}-\sigma_1\right)+\sigma_1.
\end{equation}

Under temporal concatenation, since $\frac{T}{2}>k+3$, for time $1$ to $\frac{T}{2}$, the agent will go to state $e$, stay at $e$ until time $\frac{T}{2}-2$, and go to state $f$ at time $\frac{T}{2}-1$. For time $\frac{T}{2}$ to $T-1$, the agent will have to go to $e$ again after passing through $d_1,\ldots,d_k$, stay at $e$ until time $T-2$, and then go to $f$ at time $T-1$. Recall that states $d_1$ through $d_k$ provide zero reward. The total reward for the temporal concatenation policy is
\begin{equation}
V(\mathcal I_0,\pi_{\rm TC},\mu_0)=2\left(\left(\frac{T}{2}-k\right)\left(\bar{r}-\sigma_1\right)+\sigma_1\right).
\end{equation}

Therefore, we have that the regret of temporal concatenation is given by
\begin{align}
\Delta(\mathcal I_0,\mu_0)
= & V(\mathcal I_0,\pi^*,\mu_0)-V(\mathcal I_0,\pi_{\rm TC},\mu_0)\nonumber\\
= & k\bar{r}-(k+1)\sigma_1 \nonumber\\
= & (\tau_0(\calI_0)-2)\bar{r}-(k+1)\sigma_1. \nonumber
\end{align}
By choosing $\sigma_1=\frac{\sigma}{k+1}$, we have 
\begin{equation}
\Delta(\mathcal I_0,\mu_0) =  (\tau_0(\calI_0)-2)\bar{r}-\sigma.
\end{equation}
 This completes the proof of Theorem \ref{thm:lowerb_finite_diam}. \qed

\section{Connection to the $D^*$ Diameter}
\label{sec:connection_to_d_star}
In this section we further discuss the connection between the $\epsilon$-diameter, $\tau_\epsilon$, and the diameter $D^*$ introduced in~\cite{jaksch2010near}.
As shown in~\cite{jaksch2010near}, for a time-homogeneous MDP, the value function of the optimal policy, $\pi^*$, satisfies that 
\begin{equation}
\max_{s,s'\in\mathcal S}\left(V_t^{\pi^*} (s) - V_t^{\pi^*} (s')\right) \leq \bar{r}D^*.
\end{equation}
Hence, 
\begin{equation}\label{eq:close_reward}
\sup_{\mu,\mu'\in\mathcal P}\left(V(\mathcal I, \pi^*,\mu) - V(\mathcal I, \pi^*,\mu')\right)\leq \bar{r}D^*. 
\end{equation}
Therefore, when both the reward function and the transition function are time-homogeneous, a small diameter $D^*$ guarantees that the total expected rewards of different initial distributions are close. 
However, for the more general scenarios where either the reward function or the transition function is not time-homogeneous, a small diameter $D^*$ is not sufficient for this to hold. 
In the remainder of this section, we present examples with either a time-inhomogeneous reward function, or a time-inhomogeneous transition function.
For each instance, we show that there exist two initial distributions such that the difference of total expected reward is large although the diameter $D^*$ is small. Hence, Lemma \ref{lem:upper_bound} is more general than \eqref{eq:close_reward} and the $\epsilon$-diameter more precisely characterizes the communicating property of an MDP than the diameter $D^*$.
\subsection{Time-inhomogeneous Reward Function}

In this subsection, we consider the case where the transition function is time-homogeneous but the reward function is not. We construct an example with a small $D^*$ and show that there exist two initial distributions such that the difference in total expected reward grows in $T$.

\begin{figure}[h]
\centering
\begin{tikzpicture}[->,>=stealth',auto,node distance=4cm,
  thick,main node/.style={circle,draw,font=\sffamily}]
  \node[draw=none] (0) {};
  \node[below = 0.2cm of 0]{rewards:};
  \node[main node] (1) [right of = 0] {$c$};
  \node[below = 0.05cm of 1]{$\bar{r}\mathbb I (t\in\{0,2,4,\ldots\})$};
  \node[main node] (2) [right of=1] {$d$};
   \node[below = 0.05cm of 2]{$\bar{r}\mathbb I (t\in\{1,3,5,\ldots\})$};

  \path[every node/.style={font=\sffamily\small}]
    (1) edge node [right] {} (2)
 (2) edge[bend right] node [left] {} (1);
 
\end{tikzpicture}
\caption{An MDP instance with time-inhomogeneous reward function.}
\label{fig:example_d_star_1}
\end{figure}
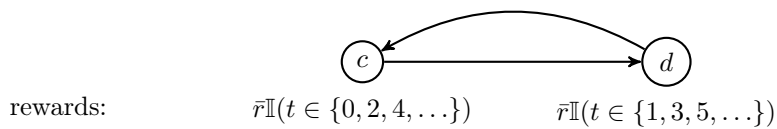

Consider the MDP instance in Figure \ref{fig:example_d_star_1}, where there are two states $c$ and $d$. The transition is deterministic: starting at state $c$, the agent can only go to $d$; starting at $d$, the agent can only go to $c$. The reward function, $R_t$, varies in time. In particular, at state $c$, the reward is $\bar{r}$ when $t$ is an even number, and $0$ when $t$ is odd; at state $d$, the reward is $\bar{r}$ when $t$ is an odd number, and $0$ when $t$ is even.

It is easy to verify that the $D^*$ diameter is small for this instance. In particular, we can go from one state to another with exactly one step. Hence, $D^* = 1$. However, we show that the initial state distribution can significantly impact the total expected reward despite the small $D^*$.
Suppose an agent starts at state $c$ when $t=0$, then the state at time $t=0,1,2,\ldots$ is deterministically $c,d,c,\ldots$, generating a reward $\bar{r}$ at each time. 
However, if the agent starts at state $d$ when $t=0$, the reward at each step is always $0$.
Therefore, 
\begin{equation}
V(\mathcal I, \pi^*,\delta_{c})- V(\mathcal I, \pi^*,\delta_{d}) = \bar{r}T,
\end{equation}
where $\delta_c$ and $\delta_d$ are the point measures at states $c$ and $d$, respectively.
Note that $\bar{r}T$ is the largest possible difference in total expected reward, which is attained in this example although $D^*$ is small.
\subsection{Time-inhomogeneous Transition Function}
\begin{figure}[h]
\centering
\includegraphics[scale = 0.8]{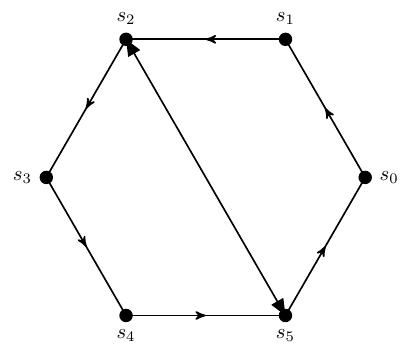}
\caption{An MDP instance with time-inhomogeneous transition function.}
\label{fig:example_d_star_2}
\end{figure}
Now we consider the case that the reward function is time-homogeneous but the transition function is not. Consider the instance in Figure \ref{fig:example_d_star_2} with states $\mathcal S = \{s_0,s_1,s_2,s_3,s_4,s_5\}$. The reward function is deterministic and depends only on the state. The reward is $\bar{r}$ for state $s_2$, and $0$ for the other five states.
The state transition is also deterministic but varies in time, which is elaborated as follows:
\begin{enumerate}
\item
Starting from state $s_0,s_1,s_3,s_4$, the system will deterministically transition to $s_1,s_2,s_4,s_5$, respectively. 
\item
Starting from state $s_2$, when $t$ is an even number, the agent can choose to go to $s_3$ or $s_5$; when $t$ is an odd number, the agent can only go to state $s_3$.
\item 
Starting from state $s_5$, the agent can choose to go to $s_0$ or $s_2$ if $t$ is odd; the system can only go to $s_0$ if $t$ is even.
\end{enumerate}

Intuitively, there is a one-way ``bridge'' between states $s_2$ and $s_5$ that shifts direction at each time step.
It is easy to see that this instance has a small $D^*$ diameter with $D^*\leq 5$.

Suppose the agent starts from state $s_2$ at $t=0$. Then the optimal policy is to go to $s_5$, then keep returning between $s_2$ and $s_5$. Note that this is feasible because whenever the agent arrives in states $s_2$ and $s_5$, the bridge is always in the proper direction such that the agent can go through it. The total expected reward can be expanded as 
\begin{equation}
V(\mathcal I,\pi^*,\delta_{s_2}) = \frac{\bar{r}T}{2}.
\end{equation}

However, if the system starts from state $s_4$, whenever the agent arrives in states $s_2$ and $s_5$, the bridge is always in the opposite direction such that the agent can never use it. Then the optimal policy is to go to $s_5,s_0,s_1,\ldots$, with 
\begin{equation}
V(\mathcal I,\pi^*,\delta_{s_4}) = \frac{\bar{r}T}{6}.
\end{equation}

Hence, the difference in total expected reward with initial distributions $\delta_{s_2}$ and $\delta_{s_5}$ is 
\begin{equation}
V(\mathcal I,\pi^*,\delta_{s_2}) - V(\mathcal I,\pi^*,\delta_{s_5}) = \frac{\bar{r}T}{3},
\end{equation}
which grows linearly in $T$ despite the small $D^*$ diameter.

\section{Proofs of Additional Theoretical Results}
\subsection{{Proof of Lemma \ref{lemma:connection of diameters}}}\label{app:lemma:connection of diameters}

Recall that $G\in\mathcal G_{\rm csl}$ is strongly connected with classical graph diameter $d_c(G)<\infty$.
For vertices $v,v'\in\mathcal V$, define ${d_G}(v,v')$ as the distance from vertex $v$ to vertex $v'$. The classical diameter satisfies $d_c(G)=\max_{v,v'\in\mathcal V}{d}_G(v,v')$. 

(1) First, we show that $\tau_0(\mathcal I_G)\geq d_c(G)$. 
By the definition of the classical diameter, there exist vertices $v_i,v_j\in\mathcal V$ with distance ${d}_G(v_i,v_j)=d_c(G)$.
Consider the point measure on state $i$, $\nu$ with $\nu(i)=1$, and the point measure on state $j$, $\nu'$ with $\nu'(j)=1$. Starting from distribution $\nu$, it takes at least $d_c(G)$ steps to achieve distribution $\nu'$, which implies that $\tau_0(\mathcal I_G)\geq d_c(G)$. 

(2) Now, we prove that $\tau_0(\mathcal I_G)\leq 2d_c(G)$ analogously to the proof of Theorem \ref{thm:lowerb_finite_diam}. Recall that $G\in\mathcal G_{\rm csl}$ has at least one self-loop. Without loss of generality, let $e_{11}\in\mathcal E$. For any pair of distributions $\nu,\nu'$, we show that starting from the initial distribution $\nu$, the system can reach $\nu'$ after $2d_c(G)$ steps. In particular, we consider the following policy:

Stage 1: If the initial state is $1$, stay at state $1$ for $d_c(G)$ steps until time $t=d_c(G)$; otherwise, if the initial state is $i$, $i\neq 1$, first go to state $1$ using ${d}_G(i,1)$ steps, and stay at state $1$ for $(d_c(G)-{d}_G(i,1))$ steps until time $t=d_c(G)$. This stage requires $d_c(G)$ steps.

Stage 2: Starting from state $1$, we reach the distribution $\nu'$ in another $d_c(G)$ steps. Analogous to Stage 2 in the proof of Theorem \ref{thm:lowerb_finite_diam}, for any state $i\in\mathcal S$, we can take a $d_c(G)$-path to reach state $i$ starting from state $1$: first stay at state $1$ for $(d_c(G)-{d}_G(1,i))$ steps, then go to state $i$ using another ${d}_G(1,i)$ steps. In order to reach the distribution $\nu'$ in exactly $d_c(G)$ steps starting from state $1$, with probability $\nu'(j)$ we take a $d_c(G)$-path to state $j$, for $j\in\mathcal S$. This stage requires $d_c(G)$ steps. 

Hence, starting from an arbitrary distribution $\nu$, the policy introduced above achieves any distribution $\nu'$ in $2d_c(G)$ steps, which implies that the $0$-diameter is at most $2d_c(G)$, i.e. $\tau_0(\calI_G)\leq 2d_c(G)$.

Combining (1) and (2) completes the proof of this lemma. \qed

\subsection{Proof of Lemma \ref{lemma:xisgt diameter}}\label{app:lemma:xisgt diameter}
Fix $\epsilon\in\big [1-(1-\xi)^{\tau_0(\mathcal I_G)},1\big )$.
It suffices to show that for any pair of distributions over the state space $\nu,\nu'\in\mathcal P$, $\nu'$ is $\epsilon$-reachable from $\nu$ in $\tau_0(\mathcal I_G)$ steps.

Recall that when $\xi=0$, the $\xi$-SGT instance becomes a DGT instance with $0$-diameter $\tau_0(\mathcal I_G)$. Hence, for any pair of distributions $\nu,\nu'$, there exists a policy $\pi_{\nu,\nu'}$ such that the state distribution transitions from $\nu$ to $\nu'$ precisely after $\tau_0(\mathcal I_G)$ steps. 

We can formulate the following policy, $\pi_{\nu,\nu'}^\xi$, for the $\xi$-SGT instance: At each time step $t$, the agent makes the decision using the same policy as in the corresponding DGT instance, pretending $\xi=0$, i.e.,
\begin{equation}
\pi_{\nu,\nu'}^\xi\left(t,S_t,Y^P\right) = \pi_{\mu_t,\nu'}\left(t,S_t,Y^P\right),
\end{equation}  
where $\mu_t$ is the state distribution at time $t$. In other words, at time $t$, the policy $\pi^\xi_{\nu,\nu'}$ makes the same decision as would $\pi_{\mu_t,\nu'}$. 

When applying the policy $\pi_{\nu,\nu'}^\xi$ in the $\xi$-SGT instance, $\mathcal I_G^\xi$, 
during the first $\tau_0(\mathcal I_G)$ steps, let $E$ be the event that the random perturbation never occurs, and $\bar{E}$ the event that the perturbation occurs in at least one step. Recall that whether the perturbation occurs in each step is independent, the probability of the event $E$ is at least $(1-\xi)^{\tau_0(\mathcal I_G)}$: 
\begin{equation}
\mathbb P(E)\geq (1-\xi)^{\tau_0(I_G)}. 
\end{equation}
Therefore, with probability at least $(1-\xi)^{\tau_0(\mathcal I_G)}$, the state distribution becomes exactly $\nu'$ after $\tau_0(\mathcal I_G)$ steps under $\pi_{\nu,\nu'}^\xi$. Otherwise, the stochastic perturbation takes place in at least one step before the state distribution reaches $\nu'$. Starting from $\nu$ at time $t=0$, the distribution at time $t=\tau_0(\mathcal I_G)$ under $\pi_{\nu,\nu'}^\xi$, {denoted by $\tilde{\mu}_{\tau_0(\mathcal I_G)}$}, can be expanded as
\begin{align}
 & \tilde{\mu}_{\tau_0(\mathcal I_G)}(s) \nonumber\\
= & \mathbb P_{S_0\sim \nu}\left(S_{\tau_0(\mathcal I_G)}=s\Big| E\right)\mathbb P(E)+ \mathbb P_{S_0\sim \nu}\left(S_{\tau_0(\mathcal I_G)}=s\Big| \bar{E}\right)\mathbb P(\bar{E})\nonumber\\
= & \nu'(s)\mathbb P(E) + \mathbb P_{S_0\sim \nu}\left(S_{\tau_0(\mathcal I_G)}=s\Big| \bar{E}\right)\mathbb P(\bar{E}),\label{pf:lem5.4-1}
\end{align}
for $s\in\mathcal S$, where Eq.~\eqref{pf:lem5.4-1} follows from the fact that the distribution of $S_{\tau_0(\mathcal I_G)}$ conditioned on the event $E$, i.e. conditioned on the event that no perturbation occurs, is exactly $\nu'$. The total variation between $\tilde{\mu}_{\tau_0(\mathcal I_G)}$ and the target distribution $\nu'$ satisfies 
\begin{align}
    & \delta_{\rm TV}\left(\tilde{\mu}_{\tau_0(\mathcal I_G)},\nu'\right) \nonumber\\
 =&\frac{1}{2}\sum_{s\in\mathcal S}\left\vert\tilde{\mu}_{\tau_0(\mathcal I_G)}(s) - \nu'(s)\right\vert\nonumber\\
 =&\frac{1}{2}\sum_{s\in\mathcal S}\left\vert -\nu'(s)(1-\mathbb P(E))+\mathbb P_{S_0\sim \nu}\left(S_{\tau_0(\mathcal I_G)}=s\Big| \bar{E}\right)\mathbb P(\bar{E})\right\vert\nonumber\\
= & (1-\mathbb P(E))\cdot\frac{1}{2}\sum_{s\in\mathcal S}\left\vert \mathbb P_{S_0\sim \nu}\left(S_{\tau_0(\mathcal I_G)}=s\Big| \bar{E}\right)-\nu'(s)\right\vert\label{pf:lem5.4-2}\\
= & (1-\mathbb P(E))\cdot \delta_{\rm TV} \left(\mathbb P_{S_0\sim \nu}\left(S_{\tau_0(\mathcal I_G)}=\cdot\Big| \bar{E}\right),\nu'\right)\nonumber\\
\leq & 1-(1-\xi)^{\tau_0(\mathcal I_G)},\label{pf:lem5.4-3}
\end{align}
where Eq.~\eqref{pf:lem5.4-2} follows from $\mathbb P(\bar{E})=1-\mathbb P(E)$, Eq.~\eqref{pf:lem5.4-3} from $\mathbb P(E)\geq (1-\xi)^{\tau_0(\mathcal I_0)}$ and $\delta_{\rm TV}(\cdot,\cdot)\leq 1$. Hence, for $\epsilon\geq 1-(1-\xi)^{\tau_0(\mathcal I_G)}$, we have shown that the policy $\pi_{\nu,\nu'}^\xi$ achieves the target distribution $\nu'$ in $\tau_0(\mathcal I_G)$ steps within a total variation distance $\epsilon$, i.e. $\tau_\epsilon\left(\mathcal I_G^\xi\right) \leq \tau_0(\mathcal I_G)$. This completes the proof. \qed

\subsection{Proof of Lemma \ref{lemma:xisgt new}}\label{app:lemma:xisgt new}
In order to prove Lemma \ref{lemma:xisgt new}, we first introduce the following lemma.
\begin{lemma}\label{lemma:xisgt diam 1}
Fix an undirected connected graph $G=(\mathcal V,\mathcal E)\in\mathcal G_{csl}$, and $\xi\in\left(0,\frac{1}{2}\right)$. Let $\mathcal I_G$ be the DGT instance characterized by $G$, and $\mathcal I_G^\xi$ be the $\xi$-SGT instance described by $(G,\xi)$. For any pair of vertices $v,v'\in\mathcal V$, denote by $d_G(v,v')$ the classical distance between $v$ and $v'$ on graph $G$. For a policy $\pi$, let $T_{v,v'}^\pi = \min\{t\geq 1:\ S_t = v' | S_0 = v\}$ be the hitting time of $v'$ starting from $v$ under $\pi$. Then there exists a policy $\pi_{v,v'}$ such that for $t>\frac{d_G(v,v')-1}{1-2\xi}$,
\begin{equation}
\mathbb P\left(T_{v,v'}^{\pi_{v,v'}} \leq t\right) \geq 1- \frac{16\xi(1-\xi)t}{((1-2\xi)t-(d_c(v,v')-1))^2}.
\end{equation}
\end{lemma}

We first prove Lemma \ref{lemma:xisgt diam 1}. Recall that $G$ is undirected and connected. For $v,v'\in\mathcal S$, consider the following policy $\pi_{v,v'}$.
Find the shortest path from $v$ to $v'$ on $G$, $p=(v_0,v_1,\ldots,v_{d-1}, v_d)$, where $d = d_G(v,v')$, $v_0 = v$, and $v_d = v'$. For each time step $t$, if the system state is on the path $p$, i.e.,~$S_t = v_i$ for some $i\in [d]$, the agent takes the action, $a_{v_i,v_{i+1}}$, to attempt to go to the next state, $v_{i+1}$, on the path $p$, and records the next state $S_{t+1}$ in the actual path $p'$; if the system state $S_t\notin p$, then the agent attempts to go back to the path $p$ by tracing back along the path $p'$ with action $a_{S_{t+1},S_t}$, and records the next state $S_{t+1}$ in $p'$. 

Let $W_i$, $i=1,2,\ldots$ be binary random variables with 
$\mathbb P(W_i=1) = 1-\mathbb P(W_i = -1)=1-\xi$. 
Going from $v$ to $v'$ under policy $\pi_{v,v'}$ corresponds to the random walk $\sum_{i=1}^t W_t$.
For each step $t$, the agent successfully achieves the intended next state with probability $1-\xi$, in which case the agent is one step closer to $v'$ under the policy $\pi_{v,v'}$, corresponding to $W_t = 1$; the agent goes elsewhere because of noise with probability $\xi$, in which case the system is one step farther from the target under $\pi_{v,v'}$, corresponding to one step back in the random walk, i.e., $W_t = -1$. Then for $t> \frac{d-1}{1-2\xi}$, we have 
\begin{align}
  \mathbb P (T_{v,v'}^{\pi_{v,v'}}\geq t) 
= & \mathbb P\left(\sum_{i = 1}^{t}W_i\leq d-1\right)\\
\leq & \frac{16\xi(1-\xi)t}{((1-2\xi)t-(d-1))^2}, \label{eq:chernoff}
\end{align}
where \eqref{eq:chernoff} follows from the Chebyshev inequality, which completes the proof of the lemma.

Now we prove Lemma \ref{lemma:xisgt new} using Lemma \ref{lemma:xisgt diam 1}.
Given $\mu,\mu'\in\mathcal P$, let $\mu' = \mu F$, where the matrix $F\in \mathbb R^{|\mathcal S|\times |\mathcal S|}$ satisfies that $F_{i,j}\in(0,1)$ for $i,j\in\mathcal S$, and $F {\bf 1} = {\bf 1}$.
Let the initial state be drawn from $\mu$, i.e., $S_0\sim\mu$.
We consider the following policy in order to achieve state distribution $\mu'$ starting from $\mu$.
If the initial state is $S_0=i\in\mathcal S$, sample a random variable $J\in\mathcal S$ with $P(J = \cdot |S_0=i) = F(i,\cdot)$. The agent then employs the policy $\pi_{i,J}$ described in the proof of Lemma \ref{lemma:xisgt diam 1} to attempt to go to state $J$, and stays at the target state $J$ after reaching it. Note that this is feasible because each state has a noiseless self-loop.
Hence, we have 
\begin{align}
\mathbb P(S_t = j| S_0=i) 
= & \mathbb P (S_t = j, J = j|S_0 = i) + \mathbb P (S_t = j, J \neq j|S_0 = i)\nonumber\\
= & \mathbb P(J=j|S_0=i) \mathbb P\left (S_t = j|J=j,S_0=i\right)+ \mathbb P (S_t = j, J \neq j|S_0 = i)\nonumber\\
= & F_{i,j}  \mathbb P\left (T_{i,j}^{\pi_{i,j}}\leq t\right)+\mathbb P (S_t = j, J \neq j|S_0 = i)\label{pf:xisgt new 1}\\
\geq & F_{i,j}  \mathbb P\left (T_{i,j}^{\pi_{i,j}}\leq t\right),\label{pf:xisgt new 2}
\end{align}
where \eqref{pf:xisgt new 1} follows from the definition of the matrix $F$ and the hitting time $T_{i,j}^{\pi_{i,j}}$, and \eqref{pf:xisgt new 2} from the fact that $\mathbb P (S_t = j, J \neq j|S_0 = i)\geq 0$.
Recall that $d_c(G) = \max_{v,v'\in\mathcal V}d_G(v,v')$ is the classical diameter of graph $G$. Then for $t> \frac{d_c(G)-1}{1-2\xi}$ and $j\in\mathcal S$, the state distribution at time $t$ satisfies 
\begin{align}
\mu_t(j) = & \mathbb P_{S_0\sim \mu} (S_t = j) \nonumber\\
= &  \sum_{i\in\mathcal S}\mathbb P(S_t = j| S_0=i) \mathbb P(S_0 = i)\nonumber\\
\geq &  \sum_{i\in\mathcal S}F_{i,j}  \mathbb P\left (T_{i,j}^{\pi_{i,j}}\leq t\right) \mu(i)\label{pf:xisgt new 3}\\
= &  \mathbb P\left (T_{i,j}^{\pi_{i,j}}\leq t\right) \sum_{i\in\mathcal S}F_{i,j}  \mu(i)\nonumber\\
= & \mathbb P\left (T_{i,j}^{\pi_{i,j}}\leq t\right)\mu'(j)\label{pf:xisgt new 4}\\
\geq & \left(1- \frac{16\xi(1-\xi)t}{((1-2\xi)t-(d_c(i,j)-1))^2}\right)\mu'(j)\label{pf:xisgt new 5}\\
\geq & \left(1- \frac{16\xi(1-\xi)t}{((1-2\xi)t-(d_c(G)-1))^2}\right)\mu'(j).\label{pf:xisgt new 6}
\end{align}
Here, \eqref{pf:xisgt new 3} follows from \eqref{pf:xisgt new 2}, \eqref{pf:xisgt new 4} from $\mu F = \mu'$, \eqref{pf:xisgt new 5} from Lemma \ref{lemma:xisgt diam 1}, and \eqref{pf:xisgt new 6} from the fact that $t>\frac{d_c(G)-1}{1-2\xi}$ and $d_c(G)\geq d_G(i,j)$.
Thus, the total variation distance between $\mu_t$ and the target distribution $\mu'$ satisfies 
\begin{align}
\delta_{\rm TV} (\mu_t,\mu') 
= &  \sum_{j\in\mathcal S: \mu'(j)\geq \mu_t(j)}\left( \mu'(j)-\mu_t(j)\right)\label{pf:xisgt new 7}\\
\leq &   \sum_{j\in\mathcal S: \mu'(j)\geq \mu_t(j)}  \frac{16\xi(1-\xi)t}{((1-2\xi)t-(d_c(G)-1))^2}\mu'(j)\label{pf:xisgt new 8}\\
\leq &  \frac{16\xi(1-\xi)t}{((1-2\xi)t-(d_c(G)-1))^2},\label{pf:xisgt new 9}
\end{align}
where \eqref{pf:xisgt new 7} follows from the definition of total variation, \eqref{pf:xisgt new 8} from \eqref{pf:xisgt new 6}, and \eqref{pf:xisgt new 9} from the fact that $\sum_{j\in\mathcal S: \mu'(j)\geq \mu_t(j)}\mu'(j)\leq 1$.
With \eqref{pf:xisgt new 9}, we have $\delta_{\rm TV} (\mu_t,\mu') \leq \epsilon$ when 
\begin{equation}
t\geq \frac{d_c(G)}{1-2\xi}+\frac{4\xi(1-\xi)}{\epsilon(1-2\xi)^2}\left(2+\sqrt{4+\frac{(1-2\xi)d_c(G)\epsilon}{\xi(1-\xi)}}\right).
\end{equation}
Therefore, 
\begin{align}
\tau_\epsilon\left(\mathcal I_G^\xi\right)\leq &\frac{d_c(G)}{1-2\xi}+\frac{4\xi(1-\xi)}{\epsilon(1-2\xi)^2}\left(2+\sqrt{4+\frac{(1-2\xi)d_c(G)\epsilon}{\xi(1-\xi)}}\right)\\
\leq & \frac{d_c(G)}{1-2\xi}+\frac{4\xi(1-\xi)}{\epsilon(1-2\xi)^2}\left(2+\sqrt{4+\frac{(1-2\xi)d_c(G)}{\xi(1-\xi)}}\right),\label{pf:eps_1}
\end{align}
with \eqref{pf:eps_1} following from $\epsilon\leq 1$,which completes the proof of Lemma \ref{lemma:xisgt new}. \qed

\subsection{Proof of Theorem \ref{thm:energy}}\label{app:thm:energy}
Denote by $\mathcal I$ an instance of dynamic energy management with storage. 
With \eqref{eq:diameter_bc} and Lemma \ref{lemma:xisgt new}, we have $\tau_\epsilon(\mathcal I)\leq \frac{\frac{B}{C}+1}{1-2\beta}+\frac{f\left(\frac{B}{C}+1,\beta\right)}{\epsilon} = \alpha+\frac{\omega}{\epsilon}$. By Theorem \ref{thm:upperbound}, we have 
\begin{equation}\label{eq:energy_bound_1}
\Delta(\mu_0,T)\leq \frac{\bar{r}}{1-\epsilon}\left(\alpha+\frac{\omega}{\epsilon}\right).
\end{equation}
Now we treat $B$, $C$, and $\beta$ as fixed parameters determined by the system and the environment, and regard $\epsilon$ as a parameter we can tune. It is easy to see that 
\begin{equation}
\frac{d}{d\epsilon}\left( \frac{\bar{r}}{1-\epsilon}\left(\alpha+\frac{\omega}{\epsilon}\right)\right) = \bar{r} \left(\frac{\alpha\epsilon^2+2\omega\epsilon-\omega}{\epsilon^2(1-\epsilon)^2}\right).
\end{equation}
Hence, the right-hand side of \eqref{eq:energy_bound_1} is minimized at $\epsilon = \frac{\sqrt{\omega^2+\omega\alpha}-\omega}{\alpha}$, which implies that 
\begin{equation}
\Delta(\mu_0,T) \leq \bar{r} \left(\sqrt{\omega+\alpha}+\sqrt{\omega}\right)^2.
\end{equation}
\qed

\section{System Specifications}\label{app:sys_spec}
For the simulations in Section \ref{subsec:runtime}, we used a desktop computer with an Intel Core i7-8700k CPU and 32GB of memory. The operating system is Windows 10 and the Matlab version is R2016b.

\end{document}